\documentclass[10pt,journal,compsoc]{IEEEtran}

\usepackage{xcolor}
\usepackage[hidelinks]{hyperref}
\usepackage{algorithm}
\usepackage{algpseudocode}
\usepackage{subcaption}
\usepackage{multirow}
\newtheorem{definition}{Definition}

\newtheorem{remark}{Remark}

\usepackage[T1]{fontenc}%

\ifCLASSINFOpdf
\else
\fi

\usepackage{amsmath}
\usepackage[cmintegrals]{newtxmath}

\usepackage{bm}
\usepackage[switch]{lineno}

\usepackage{url}

\usepackage[some]{background}
\SetBgScale{1}
\SetBgContents{\parbox{20cm}{%
  \Huge Draft:  \today\\[18cm]\rotatebox{180}{\Huge Draft:  \today}}}
\SetBgColor{red}
\SetBgAngle{270}
\SetBgOpacity{0.2}

\newcommand{\revision}[1]{{\color{black} #1}}

\begin{document}

\title{Circularity of Thermodynamical Material Networks: Indicators, Examples, and Algorithms}

\author{Federico~Zocco 
\thanks{F. Zocco is with the Centre for Intelligent Autonomous Manufacturing Systems, School of Electronics, Electrical Engineering and Computer Science and also with the Research Centre in Sustainable Energy, School of Mechanical and Aerospace Engineering, both at Queen's University Belfast, Northern Ireland, UK. Email: federico.zocco.fz@gmail.com.}%
\thanks{\emph{(Corresponding author: Federico Zocco)}}
}%

\maketitle

\begin{abstract}
\revision{The transition towards a circular economy has gained importance over the last years since the traditional linear take-make-dispose paradigm is not sustainable in the long term. Recently, thermodynamical material networks (TMNs) \cite{zocco2023thermodynamical} have been proposed as an approach to design material flows based on the idea that any supply chain can be seen as a set of thermodynamic compartments that can be added, removed, modified or connected differently. Compared to the well-established material flow analysis (MFA), TMNs leverage dynamical energy balances and ordinary differential equations along with the usual mass balances, thus tackling circular economy as a material network design problem analogous to traditional engineering design approaches (e.g., design of thermodynamic cycles, electrical and hydraulic networks) rather than as an analysis of stock-and-flow data. Hence, TMNs allow the depiction of highly dynamic material stocks and flows whose variations can occur in less than 1 minute; achieving such modelling accuracy with MFA would be more data intensive. In this paper, we first develop several circularity indicators of TMNs using a graph-based formalism. Then, we illustrate their calculation using two numerical examples for the case of fluid materials and one numerical example for the case of solid materials, for which the detailed hybrid dynamical equations and simulation outputs are provided. The paper source code is publicly available\footnotemark{}.}\footnotetext{\url{https://github.com/fedezocco/TMNcircularity-MATLAB}}
\end{abstract}


\IEEEpeerreviewmaketitle

\section{Introduction}

\revision{A ``linear'' economy extracts finite natural resources through mining and disposes of waste products in landfills. This is unsustainable in the long-term both at the extraction stage due to the finiteness of minerals and at the disposal stage due to the associated risk of pollution.} A ``circular'' economy has been proposed to address the problem by closing as much as possible the flow of materials through an increase of reuse, repair, refurbishment and recycling of goods \cite{EllenMacArthurFoundation}. As for any emerging field, defining the theoretical foundations of a circular economy is an active area of research \cite{suarez2019operational}. In particular, it is of paramount importance to clearly define the indicators of material circularity since we cannot pursue the optimization of concepts that are not clearly measurable. 

\revision{
Several proposed indicators have, at least, one of these two major limitations, which we aim to address with this paper: the first limitation is the lack of a mathematical and physical basis, which is needed for the clarity and repeatability of the measures of circularity and, eventually, for their comparisons. For example, only one algebraic equation was used by Niero \emph{et al.} \cite{niero2019coupling}, three were used by Bastianoni \emph{et al.} \cite{bastianoni2023lca}, zero were used by Pollard \emph{et al.} \cite{pollard2022developing} and by Sacco \emph{et al.} \cite{sacco2021circular}, and only two were used by Linder \emph{et al.} \cite{linder2020product}. The second major limitation is the lack of consideration of the extent to which material life-cycles are closed, which is the main focus of circularity according to the Ellen MacArthur Foundation's definition \cite{EllenMacArthurFoundation}. Examples of papers proposing circularity indicators for single life-cycle stages and/or not focused on materials are Zocco \emph{et al.} \cite{zocco2024towards}, which considers a robotic disassembly process, Toro \emph{et al.} \cite{toro2022circularity}, which considers a municipal solid waste treatment plant, Linder \emph{et al.} \cite{linder2017metric}, which looks at economic values of products, Di Maio \emph{et al.} \cite{di2015robust}, which focuses on the recycling stage, Park \emph{et al.} \cite{park2014establishing}, which addresses the share of recycled resources, and Gehin \emph{et al.} \cite{gehin2008tool}, which focuses on the stage of product design.            
}

\revision{
Static (i.e., constant) stocks and flows are the easiest model to work with, but also the least accurate. Indeed, in reality, materials move from a place to another with the consequence that the distribution of stocks and flows changes every second, that is, it is highly dynamic. Engineers and scientists typically use differential and difference equations to mathematically describe dynamic phenomena \cite{sopasakis2023control,haddad2008nonlinear}. While creating a dynamic model is more time consuming than the static counterpart, the predictions generated by the former are more realistic, and hence, more useful in simulating what-if scenarios and to answer design questions. If one or more parameters of a system are not known exactly, it is common practice to describe them with suitable probability distributions (e.g., normal, Poisson). This yields stochastic dynamic models \cite{oksendal2013stochastic,aastrom2012introduction}. In this paper, we consider static and dynamic stocks and flows in a deterministic setting. Specifically, the considered dynamics are either continuous-time or hybrid, i.e., involving both continuous- and discrete-time terms.   
}

Since to date the term ``circular economy'' has an unclear definition \cite{kirchherr2017conceptualizing}, for the sake of clarity we now state its meaning in this paper: the adjective ``circular'' means ``closed flow of material''; as a consequence, ``circular economy'' means ``an economy based on closed flows of materials'', i.e., an economy in which no input and no output material flows exist, and the expression ``measuring the network circularity'' means ``measuring to what extent the material flow in the network is closed''.

\revision{
\textbf{Our contributions:} 
\begin{enumerate}
\item{We provide a physics-based formalization and measure of material circularity by developing several graph-based indicators with their algorithms; this is achieved by leveraging thermodynamical material networks (TMNs) \cite{zocco2023thermodynamical}. Since graph theory is commonly used to mathematically describe networked systems, it can be the tool to systematically define circularity indicators.}
\item{We illustrate the applicability of the indicators using two numerical examples for the case of fluid materials (Section \ref{sec:Examples}) and a numerical example for the case of solid plastics (Section \ref{sec:caseStudy}).}
\item{The indicators are defined and calculated in the examples considering static, continuous-time dynamic, and hybrid dynamic stocks and flows.}
\item{We detail the dynamical equations of a TMN for the case of solid materials (Section \ref{sec:caseStudy}).} 
\end{enumerate}
}

The paper is organized as follows: Section \ref{sec:RelWork} discusses the related work, Section \ref{sec:InAndAlg} details the indicators and their algorithms, \revision{Section \ref{sec:Examples} covers two examples for fluid materials, while Section \ref{sec:caseStudy} considers an example for solid materials; finally, Section \ref{sec:Concl} gives the conclusions.} Throughout the paper, matrices and vectors are indicated in bold, whereas sets are indicated with calligraphic capital letters.

\section{Related Work}\label{sec:RelWork}
\subsection{Circularity Indicators} 
Over the years, indicators of circularity have been developed by different organizations and countries. For example, in 2021 the OECD published an inventory of 474 circularity indicators for different sectors such as food, water, and waste \cite{OECD}, while the European Union identified self-sufficiency of raw materials and waste generation as indicators \cite{Eurostat}, which correspond to the input and the output material flows of a linear economy, respectively. China has defined indicators that focus on different scales with the so called ``three plus one'' plan, where the ``three'' are the micro, meso, and macro levels and the ``one'' is the waste industry \cite{bleischwitz2022circular}, while the United States Environmental Protection Agency has developed tools to measure the production, the use, and the waste stream of a product \cite{US-EPA}. In this scenario, ``circular economy'' is an unclear concept incorporating different definitions \cite{kirchherr2017conceptualizing} and based on measurements of different variables such as energy flows, material flows, social impact or economic impact without well-defined mathematical and physical foundations \cite{moraga2019circular}. \revision{In this paper, we seek to add clarity to the topic by leveraging graph theory, which is a well-established mathematical formalism to model and design distributed systems such as hydraulic \cite{deuerlein2008decomposition}, electrical \cite{freitas2020stochastic}, and multiagent networks \cite{mesbahi2010graph}.}

\subsection{Graph Theory for Circular Economy} 
\revision{In Gribaudo \emph{et al.}} \cite{gribaudo2020circular}, the flow of chitin was modeled, which is in common with our use of graphs to model material flows. Graphs were used also by Singh \emph{et al.} \cite{singh2020managing} for the Indian mining industry, but not to model material flows; instead, they looked at the factors that hinder the circularity such as financial barriers, government policies, and organizational barriers. Moktadir \emph{et al.} \cite{moktadir2018drivers} modeled the relationships between four main factors that have a strong influence on the transition towards a circular economy: knowledge about circular economy, customer awareness, commitment from top management, and government legislation. The work of How \emph{et al.} \cite{how2018debottlenecking} integrates P-graphs with a sustainability index; this integration of a measure of sustainability and a material flow modeling using graphs is a key feature of our work as well. P-graphs were used also by Yeo \emph{et al.} \cite{yeo2020synthesis} to formulate and solve a combinatorial optimization problem of a biomass supply network considering the fertilizer and the electricity for recycling among other resources. The result showed that the linear economy is more profitable, whereas the circular economy could reduce the electricity import. The P-graph developed by Van \emph{et al.} \cite{van2020implementing} organizes in layers the stages of municipal solid waste processing: the input to the network is a mixed waste, thus the first layer performs the material separation, the second layer performs waste treatment and disposal, while the third layer contains the recovered products, the emissions emitted, and the emissions avoided. Gribaudo \emph{et al.} \cite{gribaudo2020circular2} proposed Petri nets to depict the end-of-life vehicles flow including stages such as demolition, recycling, scrapping, foundry, incinerators, and customers, while Hale \emph{et al.} \cite{hale2020stability} developed a dynamical social-ecological system that considers a natural resource harvested by a fixed number of consuming agents and also provide a Lyapunov-based stability proof of the system: their work is a step towards a dynamical systems theory for sustainability of social-ecological systems.

\subsection{Thermodynamical Material Networks}
TMNs were recently proposed by Zocco \emph{et al.} \cite{zocco2023thermodynamical} to provide both a generalized and a systematic approach to material flow design. The reason for such a name is: ``thermodynamical'' originates from the fact that the approach relies on compartmental dynamical thermodynamics \cite{haddad2019dynamical}; ``material'' originates from the fact that the thermodynamic compartments are connected by the flow of materials; finally, ``networks'' originates from the fact that a network is the natural result of the connection of multiple thermodynamic compartments. The generality of thermodynamics \cite{haddad2017thermodynamics,bakshi2011thermodynamics} allows material flow systems to be seen as a connection of thermodynamic compartments that can be added, removed or modified as needed to achieve a circular flow of material. This can be seen as the generalization of the design approach of hydraulic networks: the water can be replaced by any material, then the general mass and energy balances are applied to each material processing stage to describe their dynamics as hydraulic engineers do to study valves, pipes, reservoirs, and pumps. In this paper, we develop a measure of circularity for TMNs through several indicators to extend Step 1 of the four-step methodology proposed by Zocco \emph{et al.} \cite{zocco2023thermodynamical}.

\section{Indicators and Algorithms}\label{sec:InAndAlg} 
\subsection{Indicators}\label{sub:DefsProps}
We first recall some basic definitions from graph theory (see \cite{bondy1976graph}, Chapter 1), then we define several quantities related to material flow circularity.
\begin{definition}%
A \emph{directed graph} $D$, briefly called \emph{digraph}, is a graph identified by a set of \emph{vertices} (a.k.a. \emph{nodes}) $\{v_1, v_2, \dots, v_{n_v}\}$ and a set of \emph{arcs} $\{a_1, a_2, \dots, a_{n_a}\}$. \revision{An arc connects two distinct nodes and it is represented by an arrow with a direction. The node touched by the head of the arrow is the \emph{head} of the arc, while the node touched by the other end of the arc is the \emph{tail} of the arc.}   
\end{definition}
\revision{In general, more than one arc might connect two nodes.}
\begin{definition} %
A \emph{directed walk} in $D$ is a finite non-\revision{empty} sequence $W = (v_0, a_1, v_1, a_2, \dots, a_l, v_l)$ whose terms are alternately vertices and arcs such that, for $i = 1, 2, \dots, l$, the arc $a_i$ has head $v_i$ and tail $v_{i-1}$. The vertices $v_0$ and $v_l$ are the \emph{origin} and the \emph{terminus} of $W$, respectively, while $v_1, \dots, v_{l-1}$ are its \emph{internal vertices}. The integer $l$ is the \emph{length} of $W$.
\end{definition}
\begin{definition}
If the arcs $a_1, a_2, \dots, a_l$ of a directed walk $W$ are distinct, W is a \emph{directed trail}. 
\end{definition}
\begin{definition}
A directed trail is \emph{closed} if it has positive length and its origin and terminus are the same, i.e., $v_0 = v_l$. 
\end{definition}
\begin{definition}\label{def:phi}
A closed directed trail whose origin and internal vertices are distinct is a \emph{directed cycle} $\phi$. 
\end{definition}
Summarizing, a directed cycle is a directed walk $W$ in which 
\begin{enumerate}
\item{the arcs are distinct}
\item{the origin and the internal vertices are distinct}
\item{the origin and the terminus are the same}
\item{$l > 0$}.
\end{enumerate}
As we will show \revision{later in this section}, the directed cycles in a digraph are a measure of circularity if the orientation of the arcs corresponds to the flows of material within the network.

We now introduce new concepts, which are relevant to circular flow design and based on the definitions above borrowed from graph theory.

Let $c^k_{i,j}$ be the $k$-th thermodynamic compartment through which the material moves from compartment $i$ to compartment $j$. \revision{The compartment $c^k_{i,j}$ is contained inside and indicated with a \emph{control surface} in line with the design approach of thermodynamic cycles, e.g., the Rankine cycle \cite{kaminski2017introduction}.} 
\begin{definition}[Thermodynamical material network]\label{def:TMN}
\revision{A \emph{thermodynamical material network} (TMN) is a set $\mathcal{N}$ of connected thermodynamic compartments, that is, 
\begin{equation}\label{def:TMNset}
\begin{gathered}
\mathcal{N} = \left\{c^1_{1,1}, \dots, c^{k_v}_{k_v,k_v}, \dots, c^{n_v}_{n_v,n_v}, \right. \\ 
\left. c^{n_v+1}_{i_{n_v+1},j_{n_v+1}}, \dots, c^{n_v+k_a}_{i_{n_v+k_a},j_{n_v+k_a}}, \dots, c^{n_c}_{i_{n_c},j_{n_c}}\right\}, 
\end{gathered}
\end{equation}
which transport, store, use, and transform a target material. Each compartment is indicated by a \emph{control surface} and it is modeled using \emph{dynamical systems} derived from a mass balance and/or at least one of the laws of \emph{thermodynamics} \cite{haddad2019dynamical}.}
\end{definition}

Specifically, $\mathcal{N} = \mathcal{R} \cup \mathcal{T}$, where $\mathcal{R} \subseteq \mathcal{N}$ is the subset of compartments $c^k_{i,j}$ that \emph{store}, \emph{transform}, or \emph{use} the target material, while $\mathcal{T} \subset \mathcal{N}$ is the subset of compartments $c^k_{i,j}$ that \emph{move} the target material between the compartments belonging to $\mathcal{R} \subseteq \mathcal{N}$. A net $\mathcal{N}$ is associated with its weighted \emph{mass-flow digraph} $M(\mathcal{N})$, which is a weighted digraph whose nodes are the compartments $c^k_{i,j} \in \mathcal{R}$ and whose arcs are the compartments $c^k_{i,j} \in \mathcal{T}$. For node-compartments $c^k_{i,j} \in \mathcal{R}$ it holds that \revision{$i = j = k$}, whereas for arc-compartments $c^k_{i,j} \in \mathcal{T}$ it holds that $i \neq j$ because an arc moves the material from the node-compartment $c^i_{i,i}$ to the node-compartment $c^j_{j,j}$. The orientation of an arc is given by the direction of the material flow. The superscript $k$ is the identifier of each compartment. The weight assigned to a node-compartment is the mass stock $m_k$ within the corresponding compartment, whereas the weight assigned to an arc-compartment is the mass flow rate $\dot{m}_{i,j}$ from the node-compartment $c^i_{i,i}$ to the node-compartment $c^j_{j,j}$. The superscripts $k_v$ and $k_a$ in (\ref{def:TMNset}) are the $k$-th node and the $k$-th arc, respectively, while $n_c$ and $n_v$ are the total number of compartments and nodes, respectively. Since $n_a$ is the total number of arcs, it holds that $n_c = n_v + n_a$.

\begin{definition}[Material flow network] 
\revision{A \emph{material flow network} (MFN) is a TMN whose compartments are modeled using only mass balances.}
\end{definition}

\begin{definition}[Mass-flow matrix]
The \emph{mass-flow matrix} $\bm{\Gamma}(\mathcal{N})$ associated with the network $\mathcal{N}$ is
\begin{equation}\label{eq:gammaDef}
\begin{gathered}
\bm{\Gamma}(\mathcal{N}) =
\begin{bmatrix}
\gamma_{1,1} & \dots & \gamma_{1,n_v} \\
    \vdots & \ddots & \vdots \\
\gamma_{n_v,1} &  \dots & \gamma_{n_v,n_v} 
\end{bmatrix}
= \\
\begin{bmatrix}
m_1 & \dot{m}_{1,2} & \dots & \dot{m}_{1,n_v} \\
\dot{m}_{2,1} & m_2 & \dots & \dot{m}_{2,n_v} \\
\vdots & \ddots & \ddots & \vdots \\
\dot{m}_{n_v,1} & \dot{m}_{n_v,2} & \dots & m_{n_v} 
\end{bmatrix},
\end{gathered} 
\end{equation}
\revision{where $\gamma_{i,j}$ is the entry along the $i$-th row and the $j$-th column}, the entries along the main diagonal are the weights of the vertex-compartments $c^k_{i,j} \in \mathcal{R}$ (i.e., mass stocks), and the off diagonal entries are the weights of the arc-compartments $c^k_{i,j} \in \mathcal{T}$ (i.e., mass flow rates). 
\end{definition}  
Therefore, $\bm{\Gamma}(\mathcal{N})$ is a \revision{sparse} square matrix of size $n_v \times n_v$ with non-negative real elements, i.e., $\bm{\Gamma} \in \overline{\mathbb{R}}^{n_v \times n_v}_{+}$. 

From the mass conservation principle \cite{kaminski2017introduction} \revision{and in the particular case of \emph{fluid materials}}, it follows that
\begin{equation}\label{eq:relationMassFlowContinuousTime}
\frac{\text{d}}{\text{d}t}m_{k} = \sum_{i =1}^{n_v}\dot{m}_{i,k} - \sum_{j = 1}^{n_v} \dot{m}_{k,j},
\end{equation}
\revision{where $\dot{m}_{i,j} = 0$ whenever the node-compartments $c^i_{i,i}$ and $c^j_{j,j}$ are not connected by any arcs. Let $\bm{m}$ be $\bm{m} = [m_1, m_2, \dots, m_{n_v}]^\top$. Thus, equation (\ref{eq:relationMassFlowContinuousTime}) can be written in terms of the entries of $\bm{\Gamma}(\mathcal{N})$ as} 
\begin{equation}
\frac{\text{d}}{\text{d}t}\gamma_{k,k} = \sum_{\substack{i =1 \\ i \neq k}}^{n_v}\gamma_{i,k} - \sum_{\substack{j =1 \\ j \neq k}}^{n_v} \gamma_{k,j},
\end{equation}
or equivalently in vector form as
\begin{equation}\label{eq:massBalanceVectorialForm}
\begin{gathered}
\frac{\text{d}}{\text{d}t}\bm{m} = \frac{\text{d}}{\text{d}t} 
\begin{bmatrix}
m_1 \\
m_2 \\
\vdots \\
m_{n_v}
\end{bmatrix}
= \frac{\text{d}}{\text{d}t}
\begin{bmatrix}
\gamma_{1,1} \\
\gamma_{2,2} \\
\vdots \\
\gamma_{n_v,n_v} 
\end{bmatrix}  
= \\
\begin{bmatrix}
\sum\limits_{\substack{i =1 \\ i \neq 1}}^{n_v}\gamma_{i,1} - \sum\limits_{\substack{j =1 \\ j \neq 1}}^{n_v} \gamma_{1,j} \\
\sum\limits_{\substack{i =1 \\ i \neq 2}}^{n_v}\gamma_{i,2} - \sum\limits_{\substack{j =1 \\ j \neq 2}}^{n_v} \gamma_{2,j} \\
\vdots \\
\sum\limits_{\substack{i =1 \\ i \neq n_v}}^{n_v}\gamma_{i,n_v} - \sum\limits_{\substack{j =1 \\ j \neq n_v}}^{n_v} \gamma_{n_v,j}
\end{bmatrix}
.
\end{gathered}
\end{equation}  

\revision{Instead, in the particular case of \emph{solid materials}, the transportation is performed in batches, and hence, the flow is calculated as
\begin{equation}
\dot{m}_{i,j} = \frac{\overline{m}_{i,j}}{T_{i,j}}, 
\end{equation}
where $T_{i,j}$ is the transportation time for moving the batch of mass $\overline{m}_{i,j}$ from the node $i$ to the node $j$, i.e., from the node compartment $c^i_{i,i}$ to the node compartment $c^j_{j,j}$. In the case of solid materials, the continuity equation of fluids does not hold, and hence, equations (\ref{eq:relationMassFlowContinuousTime})-(\ref{eq:massBalanceVectorialForm}) are not valid.
}

\revision{Now we introduce three quantities, namely, the cycle geometric mean, the cycle harmonic mean, and the cycle arithmetic mean, to get an average flow within a cycle, and hence, to quantify the intensity of the material flow in a cycle. The cycles in a digraph contain information about the circularity in the digraph, and hence, in the network represented by the digraph.}
\begin{definition}[Cycle geometric mean]\label{def:cycleGeoMean}
The \emph{cycle geometric mean} is 
\begin{equation}
\text{GM}(\phi) = \left(\prod_{\gamma_{i,j} \in \mathcal{Y}} \gamma_{i,j}\right)^\frac{1}{l}, 
\end{equation} 
\end{definition}
where
\begin{gather}\label{eq:setY}
\mathcal{Y} = \{\gamma_{i,j}|\dot{m}_{i,j} \in \phi \}.
\end{gather}
\revision{Note that $\phi$ is a sequence of nodes and arcs (see Definition \ref{def:phi}); to keep the notation compact, an abuse of notation is used in (\ref{eq:setY}). Indeed, it is the arc associated with $\dot{m}_{i,j}$, not $\dot{m}_{i,j}$, to actually belong to $\phi$.}

\begin{definition}[Cycle harmonic mean]
The \emph{cycle harmonic mean} is 
\begin{equation}
\text{HM}(\phi) = \frac{l}{\sum\limits_{\gamma_{i,j} \in \mathcal{Y}}\frac{1}{\gamma_{i,j}}}. 
\end{equation} 
\end{definition}

\begin{definition}[Cycle arithmetic mean]
The \emph{cycle arithmetic mean} is
\begin{equation}
\text{AM}(\phi) = \frac{1}{l} \sum\limits_{\gamma_{i,j} \in \mathcal{Y}} \gamma_{i,j}. 
\end{equation}
\end{definition}

\revision{We now use the three cycle-wise metrics defined above to quantify the circularity of the flows of the entire digraph, and hence, of the whole compartmental network represented by the digraph. This yields the six indicators in Definitions \ref{def:gmac}-\ref{def:amrc}, three of which are scaled quantities because normalized by the sum of all the flows in the network, whereas the other three metrics correspond to the scaled metrics without the normalization. Hence, the latter are referred to as ``total'' indicators.}
\begin{definition}[Geometric-mean scaled circularity]\label{def:gmac}
The \emph{geometric-mean scaled circularity} $\lambda_{\text{GS}}(\bm{\Gamma}) \in [0, 1]$ of the network $\mathcal{N}$ associated with the mass-flow matrix $\bm{\Gamma}(\mathcal{N})$ is 
\begin{equation}
\lambda_{\text{GS}}(\bm{\Gamma}) = \frac{\sum\limits^{n_\phi}_{k = 1} \text{GM}(\phi_k)}{\sum\limits^{n_\phi}_{k = 1} \text{GM}(\phi_k) + \sum\limits_{\gamma_{i,j} \in \mathcal{Q}}\gamma_{i,j}},
\end{equation}
\end{definition}
where
\begin{gather}
\mathcal{Q} = \{\gamma_{i,j}|\dot{m}_{i,j} \text{ does not belong to any directed cycle}\},
\end{gather}
\revision{$\phi_k$ is the $k$-th directed cycle}, and $n_\phi$ is the number of directed cycles in $M$.

\begin{definition}[Geometric-mean total circularity]\label{def:gmrc}
The \emph{geometric-mean total circularity} $\lambda_{\text{GT}}(\bm{\Gamma})$ of the network $\mathcal{N}$ associated with the mass-flow matrix $\bm{\Gamma}(\mathcal{N})$ is 
\begin{equation}
\lambda_{\text{GT}}(\bm{\Gamma}) = \sum\limits^{n_\phi}_{k = 1} \text{GM}(\phi_k).
\end{equation}
\end{definition}

\begin{definition}[Harmonic-mean scaled circularity]\label{def:hmac}
The \emph{harmonic-mean scaled circularity} $\lambda_{\text{HS}}(\bm{\Gamma}) \in [0, 1]$ of the network $\mathcal{N}$ associated with the mass-flow matrix $\bm{\Gamma}(\mathcal{N})$ is 
\begin{equation}
\lambda_{\text{HS}}(\bm{\Gamma}) = \frac{\sum\limits^{n_\phi}_{k = 1} \text{HM}(\phi_k)}{\sum\limits^{n_\phi}_{k = 1} \text{HM}(\phi_k) + \sum\limits_{\gamma_{i,j} \in \mathcal{Q}}\gamma_{i,j}}.
\end{equation}
\end{definition}

\begin{definition}[Harmonic-mean total circularity]\label{def:hmrc}
The \emph{harmonic-mean total circularity} $\lambda_{\text{HT}}(\bm{\Gamma})$ of the network $\mathcal{N}$ associated with the mass-flow matrix $\bm{\Gamma}(\mathcal{N})$ is 
\begin{equation}
\lambda_{\text{HT}}(\bm{\Gamma}) = \sum\limits^{n_\phi}_{k = 1} \text{HM}(\phi_k).
\end{equation}
\end{definition}

\begin{definition}[Arithmetic-mean scaled circularity]\label{def:amac}
The \emph{arithmetic-mean scaled circularity} $\lambda_{\text{AS}}(\bm{\Gamma}) \in [0, 1]$ of the network $\mathcal{N}$ associated with the mass-flow matrix $\bm{\Gamma}(\mathcal{N})$ is 
\begin{equation}
\lambda_{\text{AS}}(\bm{\Gamma}) = \frac{\sum\limits^{n_\phi}_{k = 1} \text{AM}(\phi_k)}{\sum\limits^{n_\phi}_{k = 1} \text{AM}(\phi_k) + \sum\limits_{\gamma_{i,j} \in \mathcal{Q}}\gamma_{i,j}}.
\end{equation}
\end{definition}

\begin{definition}[Arithmetic-mean total circularity]\label{def:amrc}
The \emph{arithmetic-mean total circularity} $\lambda_{\text{AT}}(\bm{\Gamma})$ of the network $\mathcal{N}$ associated with the mass-flow matrix $\bm{\Gamma}(\mathcal{N})$ is 
\begin{equation}
\lambda_{\text{AT}}(\bm{\Gamma}) = \sum\limits^{n_\phi}_{k = 1} \text{AM}(\phi_k).
\end{equation}
\end{definition}

\begin{remark}\label{pro:ifQempty}
If $\mathcal{Q} = \emptyset$, then $\lambda_{\text{GS}}(\bm{\Gamma}) = \lambda_{\text{HS}}(\bm{\Gamma}) = \lambda_{\text{AS}}(\bm{\Gamma}) = 1$. 
\end{remark}

\begin{remark}
If all the flows in the cycle $\phi$ are equal to a flow $\dot{m}$, that is,
\begin{equation}
\gamma_{i,j} = \dot{m} \quad \text{for } \gamma_{i,j} \in \mathcal{Y},
\end{equation}
then $\text{GM}(\phi) = \text{HM}(\phi) = \text{AM}(\phi) = \dot{m}$.     
\end{remark}

\begin{remark}
If all the flows in cycles are equal to a \revision{generic} flow $\dot{m}$, that is, 
\begin{equation}
\gamma_{i,j} = \dot{m} \quad \text{for } \gamma_{i,j} \in \mathcal{Y}_k, \text{ with } k = \{1, \dots, n_{\phi}\}, 
\end{equation}
then $\text{GM}(\phi) = \text{HM}(\phi) = \text{AM}(\phi) = \dot{m}$, $\lambda_{\text{GT}}(\bm{\Gamma}) = \lambda_{\text{HT}}(\bm{\Gamma}) = \lambda_{\text{AT}}(\bm{\Gamma}) = n_{\phi} \dot{m}$, and
\begin{equation}
\lambda_{\text{GS}}(\bm{\Gamma}) = \lambda_{\text{HS}}(\bm{\Gamma}) = \lambda_{\text{AS}}(\bm{\Gamma}) = \frac{n_{\phi}\dot{m}}{n_{\phi}\dot{m} + \sum\limits_{\gamma_{i,j} \in \mathcal{Q}}\gamma_{i,j}}. 
\end{equation}       
Moreover, if all the flows in $\mathcal{N}$ are equal to $\dot{m}$, that is, 
\begin{equation}
\gamma_{i,j} = \dot{m} \quad \text{for } i,j \in \{1, \dots, n_v\} \text{ with } i \neq j, 
\end{equation}
then $\text{GM}(\phi) = \text{HM}(\phi) = \text{AM}(\phi) = \dot{m}$, $\lambda_{\text{GT}}(\bm{\Gamma}) = \lambda_{\text{HT}}(\bm{\Gamma}) = \lambda_{\text{AT}}(\bm{\Gamma}) = n_{\phi} \dot{m}$, and
\begin{equation}
\lambda_{\text{GS}}(\bm{\Gamma}) = \lambda_{\text{HS}}(\bm{\Gamma}) = \lambda_{\text{AS}}(\bm{\Gamma}) = \frac{n_{\phi}}{n_{\phi} + |\mathcal{Q}|}, 
\end{equation}      
where $|\cdot|$ is the cardinality operator.
\end{remark}

\revision{The indicators above quantify the flow in directed cycles. In contrast, the following indicators capture several properties of the topology of the network, namely, the connectivity (Definition \ref{def:avConn}), the cyclicity (Definition \ref{def:cyclicity}), the sharing of flows between multiple cycles (Definition \ref{def:mfshare}), and the directionality of the flows (Definition \ref{def:direction}). In particular, the latter is highly sensitive to the numbering assigned to the nodes.}   
\begin{definition}[Average connectivity]\label{def:avConn}
The \emph{average connectivity} $\lambda_{\text{C}}(\bm{\Gamma})$ of the network $\mathcal{N}$ associated with the mass-flow matrix $\bm{\Gamma}(\mathcal{N})$ is 
\begin{equation}
\lambda_{\text{C}}(\bm{\Gamma}) = \frac{1}{n_v} \sum\limits^{n_v}_{i = 1} \text{deg}_{\text{in}}(v_i) + \text{deg}_{\text{out}}(v_i),
\end{equation}
\end{definition}
where $\text{deg}_{\text{in}}(v_i)$ and $\text{deg}_{\text{out}}(v_i)$ count the arcs with head and tail the $i$-th vertex, respectively.

\begin{definition}[Cyclicity]\label{def:cyclicity}
The \emph{cyclicity} $\lambda_{\text{Y}}(\bm{\Gamma})$ of the network $\mathcal{N}$ associated with the mass-flow matrix $\bm{\Gamma}(\mathcal{N})$ is 
\begin{equation}
\lambda_{\text{Y}}(\bm{\Gamma}) = n_{\phi}.
\end{equation}
\end{definition}

\begin{definition}[Flow sharing]\label{def:mfshare}
The \emph{flow sharing} $\lambda_{\text{S}}(\bm{\Gamma})$ of the network $\mathcal{N}$ associated with the mass-flow matrix $\bm{\Gamma}(\mathcal{N})$ is 
\begin{equation}
\lambda_{\text{S}}(\bm{\Gamma}) = \sum\limits_{\gamma_{i,j} \in \mathcal{S}}\gamma_{i,j},
\end{equation}
\end{definition}
where 
\begin{equation}
\mathcal{S} = \{\gamma_{i,j}|\dot{m}_{i,j} \text{ is shared among 2 or more cycles}\}.
\end{equation}

\begin{definition}[Directionality]\label{def:direction}
The \emph{directionality} $\lambda_{\text{D}}(\bm{\Gamma})$ of the network $\mathcal{N}$ associated with the mass-flow matrix $\bm{\Gamma}(\mathcal{N})$ is 
\begin{equation}
\lambda_{\text{D}}(\bm{\Gamma}) = \frac{\sum\limits_{i < j}\gamma_{i,j}}{\sum\limits_{i > j}\gamma_{i,j}}.
\end{equation}
\end{definition}
\revision{The directionality measures to what extent the flow in the network has the direction given by the compartments with \emph{increasing} indices, e.g., from $c^1_{1,1}$ to $c^3_{3,3}$. If $\lambda_{\text{D}}(\bm{\Gamma}) < 1$, it means that the dominant flow has the direction given by the compartments with \emph{decreasing} indices, e.g., from $c^3_{3,3}$ to $c^1_{1,1}$. Hence, $\lambda_{\text{D}}(\bm{\Gamma})$ is highly sensitive to the indices assigned to the node-compartments.}

While the concept of ``circularity'' focuses primarily on flow properties and cycles (since, as stated previously, in this paper ``circularity'' means ``closed flows of material''), further information about the material distribution across the network can be given by the following auxiliary indicators.
\begin{definition}[Total stock]\label{def:totStock}
The \emph{total stock} $\theta_{\text{S}}(\bm{\Gamma})$ of the network $\mathcal{N}$ associated with the mass-flow matrix $\bm{\Gamma}(\mathcal{N})$ is 
\begin{equation}
\theta_{\text{S}}(\bm{\Gamma}) = \sum\limits_{i = 1}^{n_v}\gamma_{i,i}.
\end{equation}
\end{definition}

\begin{definition}[Total flow]\label{def:totFlow}
The \emph{total flow} $\theta_{\text{F}}(\bm{\Gamma})$ of the network $\mathcal{N}$ associated with the mass-flow matrix $\bm{\Gamma}(\mathcal{N})$ is 
\begin{equation}
\theta_{\text{F}}(\bm{\Gamma}) = \sum\limits_{i \neq j}\gamma_{i,j}.
\end{equation}
\end{definition}

\begin{definition}[Stock distribution]\label{def:stockDistr}
The \emph{stock distribution} $\theta_{\text{D}}(\bm{\Gamma})$ of the network $\mathcal{N}$ associated with the mass-flow matrix $\bm{\Gamma}(\mathcal{N})$ is 
\begin{equation}
\begin{gathered}
\theta_{\text{D}}(\bm{\Gamma}) = \sqrt{\frac{1}{n_v - 1} \sum_{i=1}^{n_v} (\gamma_{i,i} - \mu)^2} = \\ 
\sigma(\gamma_{1,1}, \dots, \gamma_{n_v,n_v}),
\end{gathered}
\end{equation} 
\end{definition}
where $\mu = \frac{1}{n_v} \sum_{i=1}^{n_v} \gamma_{i,i}$ and $\sigma(\cdot)$ is the \revision{sample} standard deviation.

\begin{definition}[Accumulation-depletion vector]\label{def:AccDeplVect}
The \emph{accumulation-depletion vector} $\bm{\theta}_{\text{A}}(\bm{\Gamma})$ of the network $\mathcal{N}$ associated with the mass-flow matrix $\bm{\Gamma}(\mathcal{N})$ is 
\begin{equation}
\bm{\theta}_{\text{A}}(\bm{\Gamma}) = \frac{\text{d}}{\text{d}t}\bm{m}.
\end{equation} 
\end{definition}

The circularity and auxiliary indicators are summarized in Table \ref{tab:summaryIndicators}. \revision{Their graphical representation is given in Fig. \ref{fig:graphicalIndSummary} for an exemplar compartmental digraph with $n_c$ = 11, $n_v$= 5, $n_a$ = 6, and $n_\phi$ = 2. Note that all arcs in Fig. \ref{fig:graphicalLamS} (the case of $\lambda_{\text{S}}(\bm{\Gamma})$) are black because no flows are shared by the cycles. Note also that $\bm{\theta}_{\text{A}}(\bm{\Gamma})$ is omitted from Fig. \ref{fig:graphicalIndSummary} due to space constraints.}      
\begin{table*}
\centering
\caption{Summary of the circularity and auxiliary indicators.}
  \begin{tabular}{c|ccc}
Type & Name & Formula & Full definition\\    
\hline \\
\multirow{34}{*}{Circularity} & Geometric-mean scaled circularity & $\lambda_{\text{GS}}(\bm{\Gamma}) = \frac{\sum\limits^{n_\phi}_{k = 1} \text{GM}(\phi_k)}{\sum\limits^{n_\phi}_{k = 1} \text{GM}(\phi_k) + \sum\limits_{\gamma_{i,j} \in \mathcal{Q}}\gamma_{i,j}}$ & Definition \ref{def:gmac} \\ [0.7cm]
& Geometric-mean total circularity & $\lambda_{\text{GT}}(\bm{\Gamma}) = \sum\limits^{n_\phi}_{k = 1} \text{GM}(\phi_k)$ & Definition \ref{def:gmrc} \\ [0.7cm]
& Harmonic-mean scaled circularity & $\lambda_{\text{HS}}(\bm{\Gamma}) = \frac{\sum\limits^{n_\phi}_{k = 1} \text{HM}(\phi_k)}{\sum\limits^{n_\phi}_{k = 1} \text{HM}(\phi_k) + \sum\limits_{\gamma_{i,j} \in \mathcal{Q}}\gamma_{i,j}}$ & Definition \ref{def:hmac} \\ [0.7cm]
& Harmonic-mean total circularity & $\lambda_{\text{HT}}(\bm{\Gamma}) = \sum\limits^{n_\phi}_{k = 1} \text{HM}(\phi_k)$ & Definition \ref{def:hmrc} \\ [0.7cm]  
& Arithmetic-mean scaled circularity & $\lambda_{\text{AS}}(\bm{\Gamma}) = \frac{\sum\limits^{n_\phi}_{k = 1} \text{AM}(\phi_k)}{\sum\limits^{n_\phi}_{k = 1} \text{AM}(\phi_k) + \sum\limits_{\gamma_{i,j} \in \mathcal{Q}}\gamma_{i,j}}$ & Definition \ref{def:amac} \\ [0.7cm]
& Arithmetic-mean total circularity & $\lambda_{\text{AT}}(\bm{\Gamma}) = \sum\limits^{n_\phi}_{k = 1} \text{AM}(\phi_k)$ & Definition \ref{def:amrc} \\ [0.7cm]  
& Average connectivity & $\lambda_{\text{C}}(\bm{\Gamma}) = \frac{1}{n_v} \sum\limits^{n_v}_{i = 1} \text{deg}_{\text{in}}(v_i) + \text{deg}_{\text{out}}(v_i)$ & Definition \ref{def:avConn} \\ [0.7cm]
& Cyclicity & $\lambda_{\text{Y}}(\bm{\Gamma}) = n_{\phi}$ & Definition \ref{def:cyclicity} \\ [0.7cm]
& Flow sharing & $\lambda_{\text{S}}(\bm{\Gamma}) = \sum\limits_{\gamma_{i,j} \in \mathcal{S}}\gamma_{i,j}$ & Definition \ref{def:mfshare} \\ [0.7cm]
& Directionality & $\lambda_{\text{D}}(\bm{\Gamma}) = \frac{\sum\limits_{i < j}\gamma_{i,j}}{\sum\limits_{i > j}\gamma_{i,j}}$ & Definition \ref{def:direction}\\
\hline \\
\multirow{10}{*}{Auxiliary} & Total stock & $\theta_{\text{S}}(\bm{\Gamma}) = \sum\limits_{i = 1}^{n_v}\gamma_{i,i}$ & Definition \ref{def:totStock} \\ [0.7cm]
& Total flow & $\theta_{\text{F}}(\bm{\Gamma}) = \sum\limits_{i \neq j}\gamma_{i,j}$ & Definition \ref{def:totFlow}\\ [0.7cm]
& Stock distribution & $\theta_{\text{D}}(\bm{\Gamma}) = \sigma(\gamma_{1,1}, \dots, \gamma_{n_v,n_v})$ & Definition \ref{def:stockDistr}\\ [0.7cm]
& Accumulation-depletion vector & $\bm{\theta}_{\text{A}}(\bm{\Gamma}) = \frac{\text{d}}{\text{d}t}\bm{m}$ & Definition \ref{def:AccDeplVect}\\ [0.3cm]
\hline
\end{tabular}
\label{tab:summaryIndicators}
\end{table*}   
\begin{remark}
Consider two mass-flow matrices $\bm{\Gamma}_a$ and $\bm{\Gamma}_b$ with their corresponding networks $\mathcal{N}_a$ and $\mathcal{N}_b$, respectively. Then, $\mathcal{N}_a$ and $\mathcal{N}_b$ have the same values of the indicators if $\bm{\Gamma}_a = \bm{\Gamma}_b$. In general, the vice-versa does not hold. For example, a digraph configured as a square with the arcs oriented clockwise has 4 vertices, 4 arcs, one oriented cycle, and hence, $\lambda_{\text{AS}} = 1$; if the arcs are oriented counterclockwise, $\bm{\Gamma}$ changes whereas $\lambda_{\text{AS}}$ does not.    
\end{remark}
\begin{figure*}
\begin{subfigure}{.33\textwidth}
\centering
\includegraphics[width=1\textwidth]{Figures/GraphicalExplaination\_lambdaGAHAAA}
\caption{$\lambda_{\text{GS}}$, $\lambda_{\text{HS}}$, and $\lambda_{\text{AS}}$}
\label{fig:graphicalGAHAAA}
\end{subfigure}
\begin{subfigure}{.33\textwidth}
\centering
\includegraphics[width=1\textwidth]{Figures/GraphicalExplaination\_lambdaGRHRAR}
\caption{$\lambda_{\text{GT}}$, $\lambda_{\text{HT}}$, and $\lambda_{\text{AT}}$}
\label{fig:graphicalGRHRAR}
\end{subfigure}
\begin{subfigure}{.33\textwidth}
\centering
\includegraphics[width=1\textwidth]{Figures/GraphicalExplaination\_lambdaC}
\caption{$\lambda_{\text{C}}$}
\label{fig:graphicalLamC}
\end{subfigure}
\begin{subfigure}{.33\textwidth}
\centering
\includegraphics[width=1\textwidth]{Figures/GraphicalExplaination\_lambdaY}
\caption{$\lambda_{\text{Y}}$}
\label{fig:graphicalLamY}
\end{subfigure}
\begin{subfigure}{.33\textwidth}
\centering
\includegraphics[width=1\textwidth]{Figures/GraphicalExplaination\_lambdaS}
\caption{$\lambda_{\text{S}}$}
\label{fig:graphicalLamS}
\end{subfigure}
\begin{subfigure}{.33\textwidth}
\centering
\includegraphics[width=1\textwidth]{Figures/GraphicalExplaination\_lambdaD}
\caption{$\lambda_{\text{D}}$}
\label{fig:graphicalLamD}
\end{subfigure}
\begin{subfigure}{.33\textwidth}
\centering
\includegraphics[width=1\textwidth]{Figures/GraphicalExplaination\_thetaS}
\caption{$\theta_{\text{S}}$}
\label{fig:graphicalTetS}
\end{subfigure}
\begin{subfigure}{.33\textwidth}
\centering
\includegraphics[width=1\textwidth]{Figures/GraphicalExplaination\_thetaF}
\caption{$\theta_{\text{F}}$}
\label{fig:graphicalTetF}
\end{subfigure}
\begin{subfigure}{.33\textwidth}
\centering
\includegraphics[width=1\textwidth]{Figures/GraphicalExplaination\_thetaD}
\caption{$\theta_{\text{D}}$}
\label{fig:graphicalTetD}
\end{subfigure}
\caption{\revision{Graphical representation of the indicators in Table \ref{tab:summaryIndicators} considering a compartmental digraph with $n_c$ = 11, $n_v$= 5, $n_a$ = 6, and $n_\phi$ = 2. Legend: orange for flows (i.e., arcs) and stocks (i.e., nodes) that increase the indicator; light blue for arcs and nodes that decrease the indicator; green for arcs and nodes that could increase or decrease the indicator; black for arcs and nodes having no influence on the indicator.}}
\label{fig:graphicalIndSummary}
\end{figure*}

\subsection{Algorithms}\label{sub:Algorithm}
A MATLAB implementation of the indicators is publicly available\footnotemark[\value{footnote}]. The pseudo-codes of the implementations are shown in Algorithm \ref{alg:cycleDependent} and Algorithm \ref{alg:cycleIndependent}: the former covers the indicators which are based on properties of the directed cycles in $\mathcal{N}$, whereas the latter covers the cycle-independent indicators. While all these indicators could be useful in practice, those based on cycles are the most direct measure of material circularity since they actually assess to what extent the material flow is closed. 

The input ``selector'' to both pseudo-codes is a string which corresponds to the subscript of the desired indicator, e.g., selector = ``AS'' selects $\lambda_{\text{AS}}$ (note that this selector is for demonstration purpose only as it is not currently used in the actual code). The function ``triu$(\bm{\Gamma}, b)$'' returns the upper diagonal portion of $\bm{\Gamma}$ with zeros under the diagonal selected by $b$, e.g. $b = 0$ selects the main diagonal and $b = 1$ selects the diagonal just above the main diagonal. The functions ``indegree$(M)$'' and ``outdegree$(M)$'' compute the in-degree and the out-degree of all vertices in $M$, respectively, while $\bm{\Gamma}(:,k)$ and $\bm{\Gamma}(k,:)$ select the $k$-th column and the $k$-th row of $\bm{\Gamma}$, respectively.         
    
\begin{algorithm}[h]
\caption{Computation of cycle-dependent indicators: $\lambda_{\text{GS}}$, $\lambda_{\text{GT}}$, $\lambda_{\text{HS}}$, $\lambda_{\text{HT}}$, $\lambda_{\text{AS}}$, $\lambda_{\text{AT}}$, $\lambda_{\text{Y}}$, and $\lambda_{\text{S}}$.}
\label{alg:cycleDependent}
\hspace*{\algorithmicindent} \textbf{Input:} $\bm{\Gamma}$, selector \\
\hspace*{\algorithmicindent} \textbf{Output:} $\lambda_{\text{GS}}$, $\lambda_{\text{GT}}$, $\lambda_{\text{HS}}$, $\lambda_{\text{HT}}$, $\lambda_{\text{AS}}$, $\lambda_{\text{AT}}$, $\lambda_{\text{Y}}$ or $\lambda_{\text{S}}$
\begin{algorithmic}[1]
\State $M = \text{digraph}(\bm{\Gamma})$ 
\State $\mathcal{P} = \text{find\_cycles}(M)$
\State $\mathcal{Q} = \text{arcs\_not\_in\_cycles}(M)$
\If{selector $\in$ \{``GS'', ``GT''\}}
	 \For{each cycle $\phi_k \in \mathcal{P}$}  
	 	\State $a_k = \text{GM}(\phi_k)$
	 \EndFor
	 \State $\lambda_{\text{GT}} = \sum\limits^{n_\phi}_{k = 1} a_k$
	 \State $\lambda_{\text{GS}} = \frac{\lambda_{\text{GT}}}{\lambda_{\text{GT}} + \sum\limits_{\gamma_{i,j} \in \mathcal{Q}}\gamma_{i,j}}$ \Comment{Note that $\lambda_{\text{GS}} \in [0, 1]$}   
\EndIf
\If{selector $\in$ \{``HS'', ``HT''\}}
	 \For{each cycle $\phi_k \in \mathcal{P}$}  
	 	\State $a_k = \text{HM}(\phi_k)$
	 \EndFor
	 \State $\lambda_{\text{HT}} = \sum\limits^{n_\phi}_{k = 1} a_k$
	 \State $\lambda_{\text{HS}} = \frac{\lambda_{\text{HT}}}{\lambda_{\text{HT}} + \sum\limits_{\gamma_{i,j} \in \mathcal{Q}}\gamma_{i,j}}$ \Comment{Note that $\lambda_{\text{HS}} \in [0, 1]$}   
\EndIf
\If{selector $\in$ \{``AS'', ``AT''\}}
	 \For{each cycle $\phi_k \in \mathcal{P}$}  
	 	\State $a_k = \text{AM}(\phi_k)$
	 \EndFor
	 \State $\lambda_{\text{AT}} = \sum\limits^{n_\phi}_{k = 1} a_k$
	 \State $\lambda_{\text{AS}} = \frac{\lambda_{\text{AT}}}{\lambda_{\text{AT}} + \sum\limits_{\gamma_{i,j} \in \mathcal{Q}}\gamma_{i,j}}$ \Comment{Note that $\lambda_{\text{AS}} \in [0, 1]$}  
\EndIf
\If{selector == ``Y''}
	\State $\lambda_{\text{Y}} = |\mathcal{P}|$
\EndIf
\If{selector == ``S''}
	\State $\mathcal{S} = \text{find\_weights\_of\_shared\_arcs}(\mathcal{P})$
	\State $\lambda_{\text{S}} = \text{sum\_all\_elements\_in\_set}(\mathcal{S})$  
\EndIf
\end{algorithmic}
\end{algorithm}

\begin{algorithm}[h]
\caption{Computation of cycle-independent indicators: $\lambda_{\text{D}}$, $\lambda_{\text{C}}$, $\theta_{\text{S}}$, $\theta_{\text{F}}$, $\theta_{\text{D}}$, and $\bm{\theta}_{\text{A}}$.}
\label{alg:cycleIndependent}
\hspace*{\algorithmicindent} \textbf{Input:} $\bm{\Gamma}$, selector \\
\hspace*{\algorithmicindent} \textbf{Output:} $\lambda_{\text{D}}$, $\lambda_{\text{C}}$, $\theta_{\text{S}}$, $\theta_{\text{F}}$, $\theta_{\text{D}}$ or $\bm{\theta}_{\text{A}}$
\begin{algorithmic}[1]
\State $M = \text{digraph}(\bm{\Gamma})$
\If{selector == ``lambda\_D''}
\State $\lambda_{\text{D}} = \frac{\sum \text{triu}(\bm{\Gamma},1)}{\sum \gamma_{i,j} - \sum \text{triu}(\bm{\Gamma},0)}$
\EndIf 
\If{selector == ``C''}
\State $\lambda_{\text{C}} = \frac{1}{n_v} \sum (\text{indegree}(M) + \text{outdegree}(M))$
\EndIf 
\If{selector == ``S''}
\State $\theta_{\text{S}} = \text{trace}(\bm{\Gamma})$
\EndIf 
\If{selector == ``F''}
\State $\theta_{\text{F}} = \sum \gamma_{i,j} - \text{trace}(\bm{\Gamma})$ 
\EndIf 
\If{selector == ``theta\_D''}
\State $\theta_{\text{D}} = \text{std}(\gamma_{1,1}, \dots, \gamma_{n_v,n_v})$
\EndIf 
\If{selector == ``A''}
\State $\bm{\theta}_{\text{A}}(k) = (\sum\bm{\Gamma}(:,k) - \gamma_{k,k}) - (\sum\bm{\Gamma}(k,:) - \gamma_{k,k})$
\EndIf 
\end{algorithmic}
\end{algorithm}

\section{\revision{Examples for Fluid Materials}}\label{sec:Examples}
This section covers two examples: the first example illustrates the calculation of the indicators for a given $\mathcal{N}$ and its $\bm{\Gamma}(\mathcal{N})$ \revision{in the particular case of fluid material}; as the calculation highlights the violation of the mass conservation principle, the second example shows how to fix such a violation.
\subsection{\revision{Fluid Case} and Nonphysical Dynamic Graph}
\begin{figure*}
\begin{subfigure}{.33\textwidth}
\centering
\includegraphics[width=0.6\textwidth]{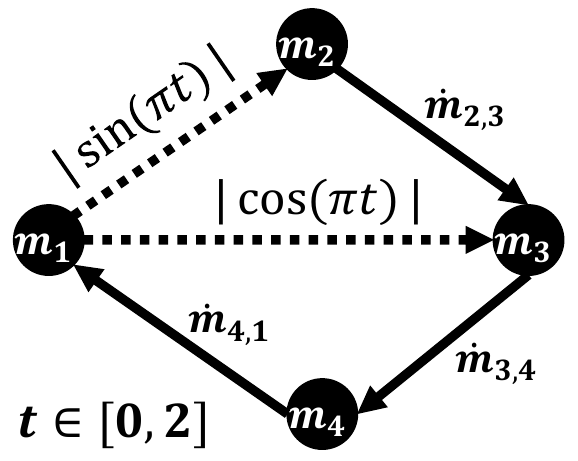}
\caption{The mass-flow digraph with its two dynamic flows (dashed arrows).}
\label{fig:Example1sub1}
\end{subfigure}
\begin{subfigure}{.33\textwidth}
\centering
\includegraphics[width=1\textwidth]{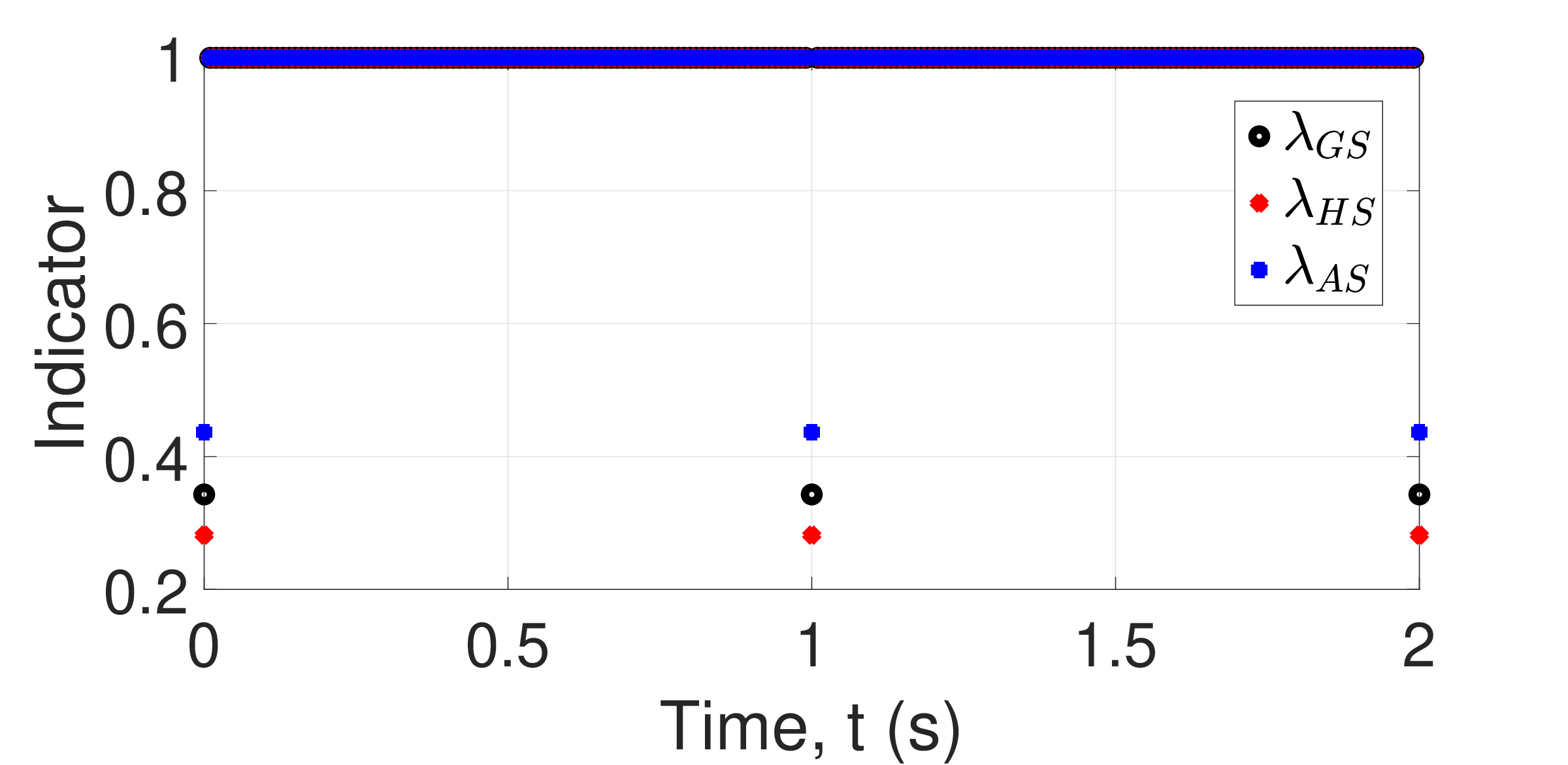}
\caption{$\lambda_{\text{GS}}(t)$, $\lambda_{\text{HS}}(t)$, and $\lambda_{\text{AS}}(t)$}
\label{fig:Example1sub2}
\end{subfigure}
\begin{subfigure}{.33\textwidth}
\centering
\includegraphics[width=1\textwidth]{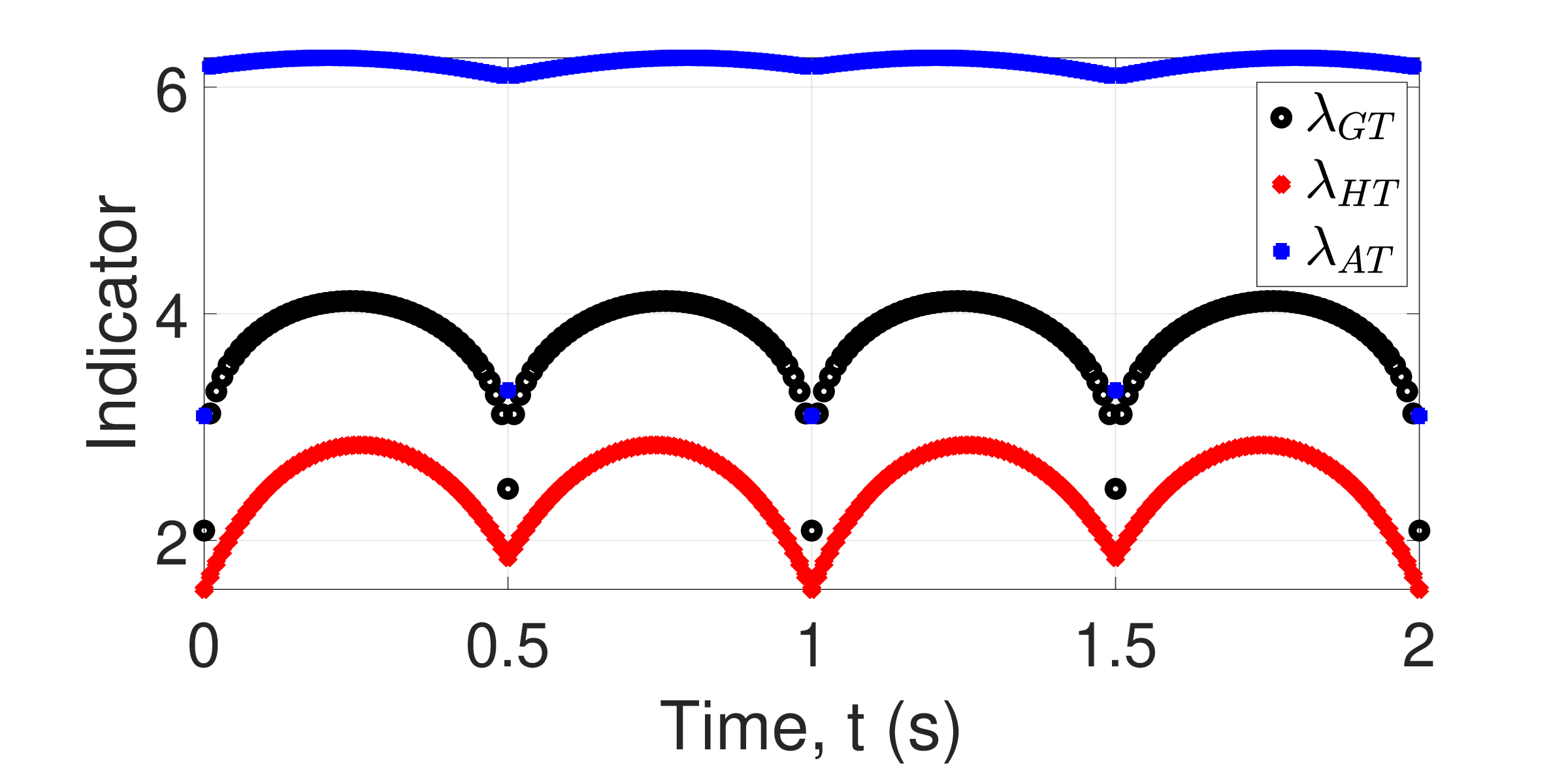}
\caption{$\lambda_{\text{GT}}(t)$, $\lambda_{\text{HT}}(t)$, and $\lambda_{\text{AT}}(t)$}
\label{fig:Example1sub3}
\end{subfigure}
\begin{subfigure}{.33\textwidth}
\centering
\includegraphics[width=1\textwidth]{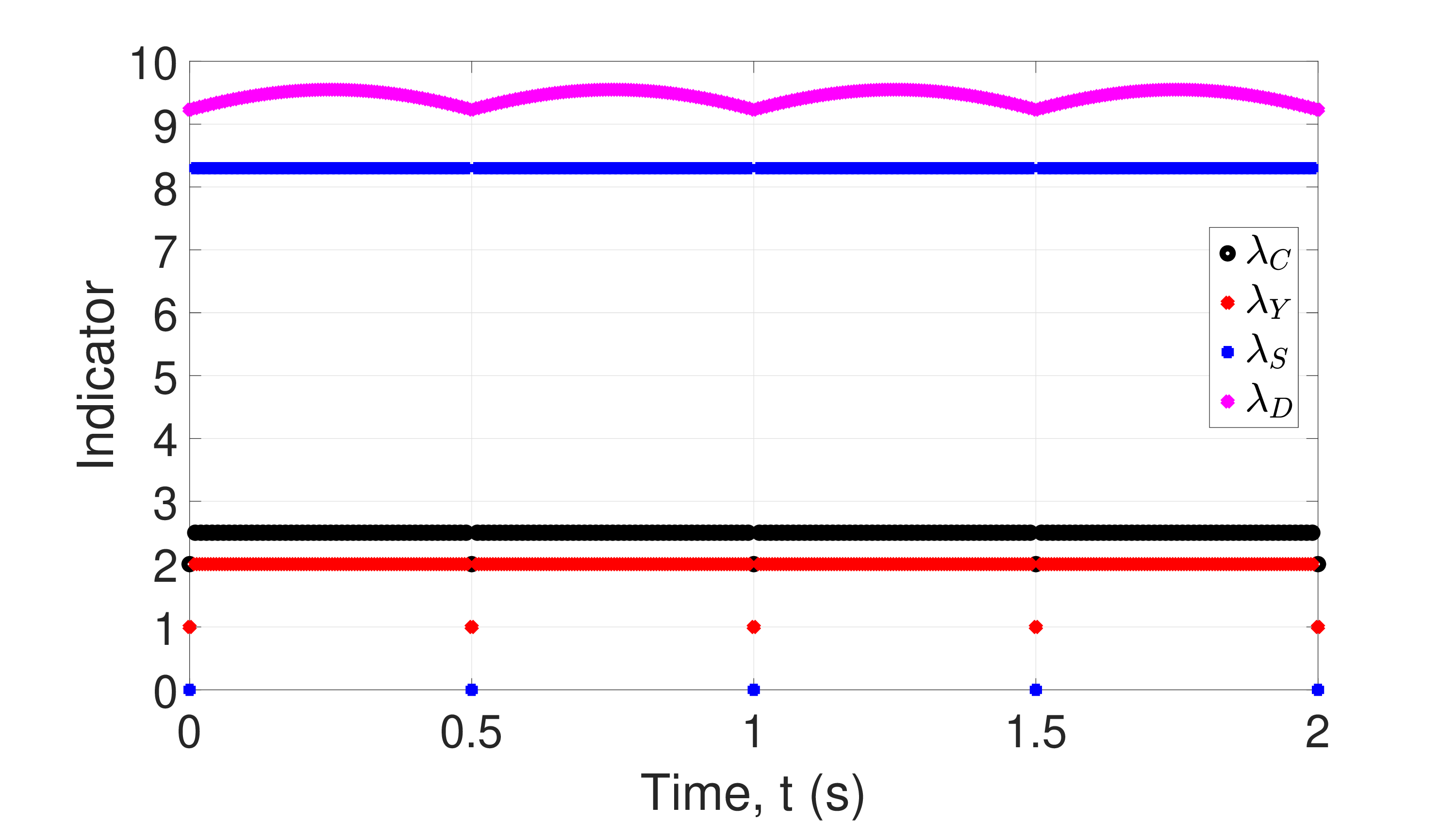}
\caption{$\lambda_{\text{C}}(t)$, $\lambda_{\text{Y}}(t)$, $\lambda_{\text{S}}(t)$, and $\lambda_{\text{D}}(t)$}
\label{fig:Example1sub4}
\end{subfigure}
\begin{subfigure}{.33\textwidth}
\centering
\includegraphics[width=1\textwidth]{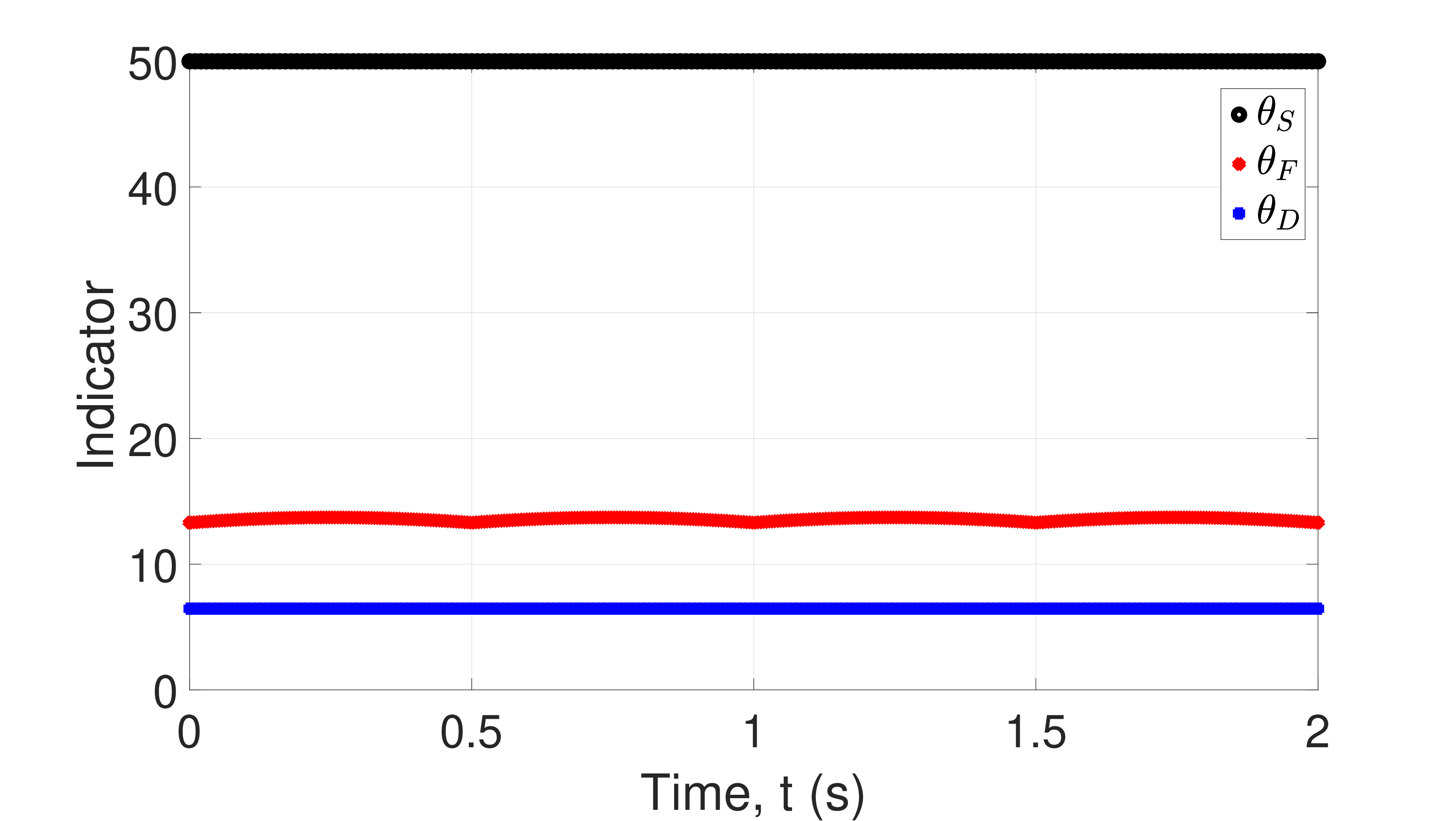}
\caption{$\theta_{\text{S}}(t)$, $\theta_{\text{F}}(t)$, and $\theta_{\text{D}}(t)$}
\label{fig:Example1sub5}
\end{subfigure}
\begin{subfigure}{.33\textwidth}
\centering
\includegraphics[width=1\textwidth]{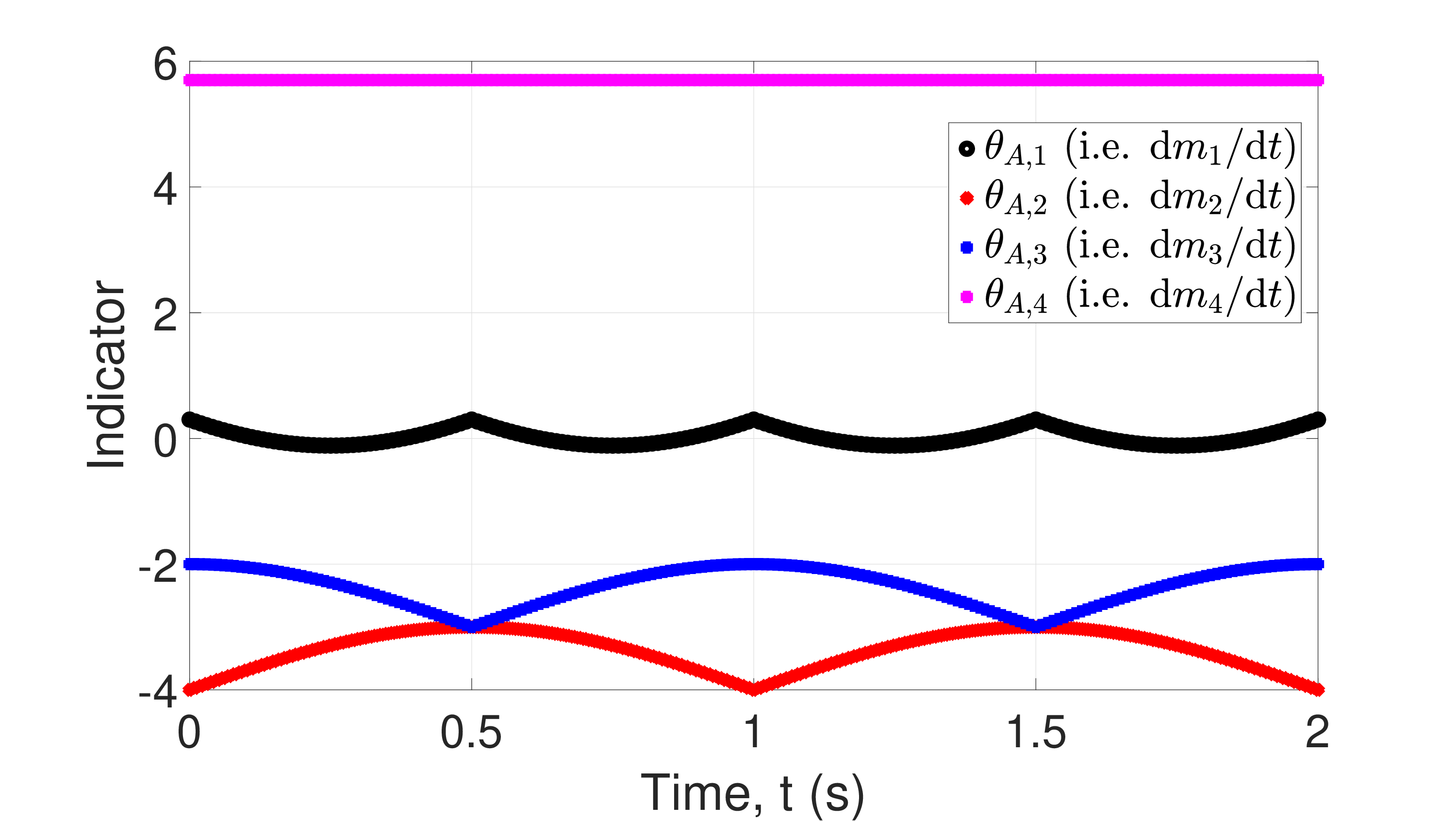}
\caption{Elements of $\bm{\theta}_{\text{A}}(t)$}
\label{fig:Example1Sub6}
\end{subfigure}
\caption{Digraph and dynamics of the indicators in Example 1 for $t \in [0, 2]$.}
\label{fig:Example1}
\end{figure*}

Assume the network 
\begin{equation}\label{eq:netExample1}
\mathcal{N} = \left\{c^1_{1,1}, c^2_{2,2}, c^3_{3,3}, c^4_{4,4}, c^5_{1,2}, c^6_{1,3}, c^7_{2,3}, c^8_{3,4}, c^9_{4,1}\right\}
\end{equation}
associated with the mass-flow matrix
\begin{equation}\label{eq:GammaExample1}
\begin{gathered}
\bm{\Gamma}(t) = 
\begin{bmatrix}
m_1 & \dot{m}_{1,2}(t) & \dot{m}_{1,3}(t) & \dot{m}_{1,4} \\
\dot{m}_{2,1} & m_2 & \dot{m}_{2,3} & \dot{m}_{2,4} \\
\dot{m}_{3,1} & \dot{m}_{3,2} & m_3 & \dot{m}_{3,4} \\
\dot{m}_{4,1} & \dot{m}_{4,2} & \dot{m}_{4,3} & m_4 \\ 
\end{bmatrix} 
=  \\
\begin{bmatrix}
10 & |\sin(\pi t)| & |\cos(\pi t)| & 0 \\
0 & 20 & 4 & 0 \\
0 & 0 & 15 & 7 \\
1.3 & 0 & 0 & 5
\end{bmatrix}
\end{gathered}
\end{equation}
with $t \in [0, 2]$. The mass-flow digraph $M(\mathcal{N})$ is shown in Fig. \ref{fig:Example1sub1}.
A topological change in $\mathcal{N}$ occurs at time $t^*$ when an $\dot{m}_{i,j}(t)$ verifies either 
\begin{equation}
\dot{m}_{i,j}(t) =
\begin{cases}
z, \quad t < t^* \\
0, \quad t = t^*  
\end{cases}
\text{or} \quad \dot{m}_{i,j}(t) =
\begin{cases}
0, \quad t < t^* \\
z, \quad t = t^*  
\end{cases},
\end{equation}
with $z > 0$. In this example, $t^* \in \{0, 0.5, 1, 1.5, 2\}$. Specifically, for $t^* = 0, 1, 2$ the net (\ref{eq:netExample1}) becomes
\begin{equation}
\mathcal{N}_1 = \left\{c^1_{1,1}, c^2_{2,2}, c^3_{3,3}, c^4_{4,4}, c^6_{1,3}, c^7_{2,3}, c^8_{3,4}, c^9_{4,1}\right\}
\end{equation}
while for $t^* = 0.5, 1.5$ it becomes
\begin{equation}
\mathcal{N}_2 = \left\{c^1_{1,1}, c^2_{2,2}, c^3_{3,3}, c^4_{4,4}, c^5_{1,2}, c^7_{2,3}, c^8_{3,4}, c^9_{4,1}\right\}.
\end{equation}  
The cases are detailed below.

\textbf{Case $t^* = 0, 1, 2$:} In this case $\dot{m}_{1,2}(t) = |\sin(\pi t)| = 0$ and $\dot{m}_{1,3}(t) = |\cos(\pi t)| = 1$, hence $\mathcal{N}_1$ has only one directed cycle
\begin{equation}\label{eq:cycle1}
\phi_1 = \left(c^1_{1,1}, c^6_{1,3}, c^3_{3,3}, c^8_{3,4}, c^4_{4,4}, c^9_{4,1}, c^1_{1,1} \right). 
\end{equation}
The indicators become:
\begin{equation}
\begin{gathered}
\lambda_{\text{GT}}(t) = \text{GM}(\phi_1) = \left(\dot{m}_{1,3}(t) \dot{m}_{3,4} \dot{m}_{4,1}\right)^{\frac{1}{3}}|_{t = 0, 1, 2} = \\ \left(9.1 |\cos(\pi t)| \right)^{\frac{1}{3}}|_{t = 0, 1, 2} = 2.1,
\end{gathered}
\end{equation}

\begin{equation}
\lambda_{\text{GS}}(t) = \frac{\lambda_{\text{GT}}(t)}{\lambda_{\text{GT}}(t) + \dot{m}_{2,3}}|_{t = 0, 1, 2} = 0.34,
\end{equation}

\begin{equation}
\begin{gathered}
\lambda_{\text{HT}}(t) = \text{HM}(\phi_1) = \frac{3}{\frac{1}{\dot{m}_{1,3}(t)} + \frac{1}{\dot{m}_{3,4}} + \frac{1}{\dot{m}_{4,1}}}|_{t = 0, 1, 2} = \\
\frac{3}{\frac{1}{|\cos(\pi t)|} + \frac{1}{7} + \frac{1}{1.3}}|_{t = 0, 1, 2} = 1.6,
\end{gathered}
\end{equation}

\begin{equation}
\lambda_{\text{HS}}(t) = \frac{\lambda_{\text{HT}}(t)}{\lambda_{\text{HT}}(t) + \dot{m}_{2,3}}|_{t = 0, 1, 2} = 0.28,
\end{equation}

\begin{equation}
\begin{gathered}
\lambda_{\text{AT}}(t) = \text{AM}(\phi_1) = \frac{1}{3} \left(\dot{m}_{1,3}(t) + \dot{m}_{3,4} + \dot{m}_{4,1}\right)|_{t = 0, 1, 2} = \\ 
\frac{1}{3}\left(|\cos(\pi t)| + 8.3\right)|_{t = 0, 1, 2} = 3.1,  
\end{gathered}
\end{equation}

\begin{equation}
\lambda_{\text{AS}}(t) = \frac{\lambda_{\text{AT}}(t)}{\lambda_{\text{AT}}(t) + \dot{m}_{2,3}}|_{t = 0, 1, 2} = 0.44, 
\end{equation}

\begin{equation}
\begin{gathered}
\lambda_{\text{C}}(t) =
\frac{1}{4} \left[\text{deg}_{\text{in}}(v_1) + \text{deg}_{\text{out}}(v_1) + \text{deg}_{\text{in}}(v_2) + \right.\\ 
\left. \text{deg}_{\text{out}}(v_2) + \text{deg}_{\text{in}}(v_3) + \text{deg}_{\text{out}}(v_3) + \text{deg}_{\text{in}}(v_4) \right. \\ 
\left. + \text{deg}_{\text{out}}(v_4)\right] = 2,
\end{gathered}
\end{equation}
$\lambda_{\text{Y}}(t) = 1$, $\lambda_{\text{S}}(t)|_{\mathcal{S} = \emptyset} = 0$,
\begin{equation}
\begin{gathered}
\lambda_{\text{D}}(t) = \frac{\dot{m}_{1,3}(t) + \dot{m}_{2,3} + \dot{m}_{3,4}}{\dot{m}_{4,1}}|_{t = 0, 1, 2} = \\ 
\frac{|\cos(\pi t)| + 11}{1.3}|_{t = 0, 1, 2} = 9.23, 
\end{gathered}
\end{equation}

\begin{equation}
\theta_{\text{S}} = m_1 + m_2 + m_3 + m_4 = 50,
\end{equation}  

\begin{equation}
\begin{gathered}
\theta_{\text{F}}(t) = \dot{m}_{1,3}(t) + \dot{m}_{2,3} + \dot{m}_{3,4} + \\ 
\dot{m}_{4,1}|_{t = 0, 1, 2} = |\cos(\pi t)| + 12.3|_{t = 0, 1, 2} = 13.3, 
\end{gathered}
\end{equation}
$\theta_{\text{D}} = 6.45$ and
\begin{equation}\label{eq:thetaAEx1Case1}
\begin{gathered}
\bm{\theta}_{\text{A}}(t) = 
\begin{bmatrix}
\dot{m}_{4,1} - |\cos(\pi t)|\\
- \dot{m}_{2,3}\\
|\cos(\pi t)| + \dot{m}_{2,3} - \dot{m}_{3,4}\\
\dot{m}_{3,4} - \dot{m}_{4,1}
\end{bmatrix} 
_{t = 0, 1, 2}=
\begin{bmatrix}
0.3 \\
-4 \\
-2 \\
5.7
\end{bmatrix}.
\end{gathered}
\end{equation}

\textbf{Case $t^* = 0.5, 1.5$:} In this case $\dot{m}_{1,3}(t) = |\cos(\pi t)| = 0$ and $\dot{m}_{1,2}(t) = |\sin(\pi t)| = 1$, hence $\mathcal{N}_2$ has only one directed cycle
\begin{equation}\label{eq:cycle2}
\phi_2 = \left(c^1_{1,1}, c^5_{1,2}, c^2_{2,2}, c^7_{2,3}, c^3_{3,3}, c^8_{3,4}, c^4_{4,4}, c^9_{4,1}, c^1_{1,1} \right).
\end{equation}
The indicators become:
\begin{equation}
\begin{gathered}
\lambda_{\text{GT}}(t) = \text{GM}(\phi_2) = \left(|\sin(\pi t)| \dot{m}_{2,3} \dot{m}_{3,4} \dot{m}_{4,1}\right)^{\frac{1}{4}}|_{t=0.5, 1.5} = \\ \left(36.4 |\sin(\pi t)| \right)^{\frac{1}{4}}|_{t=0.5, 1.5} = 2.46,  
\end{gathered}
\end{equation}
from Proposition \ref{pro:ifQempty} it follows that $\lambda_{\text{GS}}(t) = 1$, $\lambda_{\text{HS}}(t) = 1$, and $\lambda_{\text{AS}}(t) = 1$; 
\begin{equation}
\begin{gathered}
\lambda_{\text{HT}}(t) = \text{HM}(\phi_2) = \frac{4}{\frac{1}{\dot{m}_{1,2}(t)} + \frac{1}{\dot{m}_{2,3}} + \frac{1}{\dot{m}_{3,4}} + \frac{1}{\dot{m}_{4,1}}}|_{t=0.5, 1.5} = \\
\frac{4}{\frac{1}{|\sin(\pi t)|} + \frac{1}{4} + \frac{1}{7} + \frac{1}{1.3}}|_{t=0.5, 1.5} = 1.85,
\end{gathered}  
\end{equation}
\begin{equation}
\begin{gathered}
\lambda_{\text{AT}}(t) = \text{AM}(\phi_2) = \frac{1}{4}\left(\dot{m}_{1,2}(t) + \dot{m}_{2,3} + \dot{m}_{3,4} + \right.
\\\left. \dot{m}_{4,1} \right)|_{t=0.5, 1.5} = \frac{1}{4}\left(|\sin(\pi t)| + 12.3 \right)|_{t=0.5, 1.5} = 3.32,    
\end{gathered}
\end{equation}
$\lambda_{\text{C}}(t) = 2$, $\lambda_{\text{Y}}(t) = 1$, $\lambda_{\text{S}}(t)|_{\mathcal{S} = \emptyset} = 0$,
\begin{equation}
\begin{gathered}
\lambda_{\text{D}}(t) = \frac{\dot{m}_{1,2}(t) + \dot{m}_{2,3} + \dot{m}_{3,4}}{\dot{m}_{4,1}}|_{t=0.5, 1.5} = \\ 
\frac{|\sin(\pi t)| + 11}{1.3}|_{t=0.5, 1.5} = 9.23,
\end{gathered}
\end{equation} 
$\theta_{\text{S}} = 50$, 
\begin{equation}
\begin{gathered}
\theta_{\text{F}}(t) = \dot{m}_{1,2}(t) + \dot{m}_{2,3} + \dot{m}_{3,4} + \\ 
\dot{m}_{4,1}|_{t = 0.5, 1.5} = |\sin(\pi t)| + 12.3|_{t = 0.5, 1.5} = 13.3,
\end{gathered}
\end{equation} 
$\theta_{\text{D}} = 6.45$ and
\begin{equation}\label{eq:thetaAEx1Case2}
\bm{\theta}_{\text{A}}(t) = 
\begin{bmatrix}
\dot{m}_{4,1} - |\sin(\pi t)| \\
|\sin(\pi t)| - \dot{m}_{2,3} \\
\dot{m}_{2,3} - \dot{m}_{3,4} \\
\dot{m}_{3,4} - \dot{m}_{4,1}
\end{bmatrix}
_{t = 0.5, 1.5} =
\begin{bmatrix}
0.3 \\
-3 \\
-3 \\
5.7
\end{bmatrix}
.
\end{equation}

\textbf{Case $t \neq t^*$:} In this case $\dot{m}_{1,2}(t) = |\sin(\pi t)| \neq 0$ and $\dot{m}_{1,3}(t) = |\cos(\pi t)| \neq 0$, hence $\mathcal{N}$ has two cycles: $\phi_1$ as in (\ref{eq:cycle1}) and $\phi_2$ as in (\ref{eq:cycle2}). The indicators become:
\begin{equation}
\begin{gathered}
\lambda_{\text{GT}}(t) = \text{GM}(\phi_1) + \text{GM}(\phi_2) = \\
\left(9.1 |\cos(\pi t)| \right)^{\frac{1}{3}} + \left(36.4 |\sin(\pi t)| \right)^{\frac{1}{4}} 
\end{gathered}
\end{equation}
while from Proposition \ref{pro:ifQempty} it follows that $\lambda_{\text{GS}}(t) = 1$, $\lambda_{\text{HS}}(t) = 1$, and $\lambda_{\text{AS}}(t) = 1$;
\begin{equation}
\begin{gathered}
\lambda_{\text{HT}}(t) = \text{HM}(\phi_1) + \text{HM}(\phi_2) = \\
\frac{3}{\frac{1}{|\cos(\pi t)|} + \frac{1}{7} + \frac{1}{1.3}} + \frac{4}{\frac{1}{|\sin(\pi t)|} + \frac{1}{4} + \frac{1}{7} + \frac{1}{1.3}},
\end{gathered} 
\end{equation}
\begin{equation}
\begin{gathered}
\lambda_{\text{AT}}(t) = \text{AM}(\phi_1) + \text{AM}(\phi_2) = \\
\frac{1}{3}\left(|\cos(\pi t)| + 8.3\right) + \frac{1}{4}\left(|\sin(\pi t)| + 12.3 \right), 
\end{gathered}
\end{equation}
$\lambda_{\text{C}}(t) = 2.5$, $\lambda_{\text{Y}}(t) = 2$, $\lambda_{\text{S}}(t) = \dot{m}_{3,4} + \dot{m}_{4,1} = 8.3$,
\begin{equation}
\lambda_{\text{D}}(t) = 8.46 + \frac{|\sin(\pi t)| + |\cos(\pi t)|}{1.3},  
\end{equation}
$\theta_{\text{S}} = 50$, 
\begin{equation}
\theta_{\text{F}}(t) = 12.3 + |\sin(\pi t)| + |\cos(\pi t)|,
\end{equation}
$\theta_{\text{D}} = 6.45$ and
\begin{equation}\label{eq:thetaAEx1Case3}
\bm{\theta}_{\text{A}}(t) =
\begin{bmatrix}
1.3 - |\sin(\pi t)| - |\cos(\pi t)| \\
|\sin(\pi t)| - 4 \\
|\cos(\pi t)| - 3 \\
5.7
\end{bmatrix}.
\end{equation}
Now, note that the mass stocks $m_i|_{i=1,2,3,4}$ defined initially in equation (\ref{eq:GammaExample1}) are constant, which requires that  
\begin{equation}\label{eq:Ex1Violated}
\bm{\theta}_{\text{A}}(t) = \frac{\text{d}}{\text{d}t}\bm{m}(t) = \bm{0}, \quad t \in [0, 2],
\end{equation}
where $\bm{0} \in \mathbb{R}^{n_v}$ is a vector of zeros. However, equations (\ref{eq:thetaAEx1Case1}), (\ref{eq:thetaAEx1Case2}), (\ref{eq:thetaAEx1Case3}), and Fig. \ref{fig:Example1Sub6} show that (\ref{eq:Ex1Violated}) is violated. The next example solves this nonphysical behavior caused by the violation of the mass conservation principle.

\subsection{\revision{Fluid Case} and Imposition of Mass Balance}   
The violation of the mass balance in the previous example can be solved by imposing equation (\ref{eq:massBalanceVectorialForm}) and updating the entries along the main diagonal of $\bm{\Gamma}(t)$. In detail, from (\ref{eq:massBalanceVectorialForm}) it follows that
\begin{equation}\label{eq:constraintOfStocks}
\begin{gathered}
\bm{m}(t) = 
\begin{bmatrix}
m_1(t) \\
\vdots \\
m_{n_v}(t) 
\end{bmatrix} 
= 
\begin{bmatrix}
\int \sum\limits_{\substack{i =1 \\ i \neq 1}}^{n_v}\gamma_{i,1} - \sum\limits_{\substack{j =1 \\ j \neq 1}}^{n_v} \gamma_{1,j} \text{d}t\\
\int \sum\limits_{\substack{i =1 \\ i \neq 2}}^{n_v}\gamma_{i,2} - \sum\limits_{\substack{j =1 \\ j \neq 2}}^{n_v} \gamma_{2,j} \text{d}t\\
\vdots \\
\int \sum\limits_{\substack{i =1 \\ i \neq n_v}}^{n_v}\gamma_{i,n_v} - \sum\limits_{\substack{j =1 \\ j \neq n_v}}^{n_v} \gamma_{n_v,j} \text{d}t
\end{bmatrix}.
\end{gathered}
\end{equation}
The previous example is now updated including the constraint (\ref{eq:constraintOfStocks}) and the initial conditions $\bm{m}(0) = [m^0_1, m^0_2, m^0_3, m^0_4]^\top$. The following nomenclature will be adopted: $\bm{m}^{\overline{t}} = \bm{m}(\overline{t})$, $\bm{m}^{\overline{t}_a,\overline{t}_b}(t) = \bm{m}(t)|_{\overline{t}_a < t < \overline{t}_b}$, while $\bm{m}_-^{\overline{t}_a,\overline{t}_b}(t)$ and $\bm{m}_+^{\overline{t}_a,\overline{t}_b}(t)$ are the left-limit and the right-limit of $\bm{m}^{\overline{t}_a,\overline{t}_b}(t)$, respectively.

\textbf{Case $t = 0$:} 
$\bm{m}^0 = [m^0_1, m^0_2, m^0_3, m^0_4]^\top$.

\textbf{Case $0 < t < 0.5$:}
The initial condition for this interval is given by the previous interval as $\bm{m}_-^{0,0.5}(t) = \bm{m}^0$. Therefore, using (\ref{eq:thetaAEx1Case3}) and that $\sin(\pi t) > 0, \cos(\pi t) > 0$, we have that 
\begin{equation}
\bm{m}^{0,0.5}(t) = \int \bm{\theta}_{\text{A}}(t) \text{d}t  = \\
\begin{bmatrix}
1.3 t + \frac{\cos(\pi t)}{\pi} - \frac{\sin(\pi t)}{\pi} \\
- \frac{\cos(\pi t)}{\pi} - 4t \\
\frac{\sin(\pi t)}{\pi} - 3t \\
5.7 t 
\end{bmatrix}
+
\begin{bmatrix}
k^a_1 \\
k^a_2 \\
k^a_3 \\
k^a_4
\end{bmatrix},
\end{equation}
where $\bm{k}^a = [k^a_1, k^a_2, k^a_3, k^a_4]^\top = [m^0_1 - \frac{1}{\pi}, m^0_2 + \frac{1}{\pi}, m^0_3, m^0_4]$.

\textbf{Case $t = 0.5$:}
Using (\ref{eq:thetaAEx1Case2}) and that $\sin(\frac{\pi}{2}) > 0$, it follows that
\begin{equation}
\bm{m}^{0.5} = 
\begin{bmatrix}
\dot{m}_{4,1}t + \frac{\cos(\pi t)}{\pi} \\
-\frac{\cos(\pi t)}{\pi} - \dot{m}_{2,3}t \\
\dot{m}_{2,3}t - \dot{m}_{3,4}t \\
\dot{m}_{3,4}t - \dot{m}_{4,1}t
\end{bmatrix}_{t=0.5}
+
\begin{bmatrix}
k_1^b \\
k_2^b \\
k_3^b \\
k_4^b
\end{bmatrix}
=
\begin{bmatrix}
m_1^0 + 0.65 - \frac{2}{\pi} \\
m_2^0 - 2 + \frac{1}{\pi} \\
m_3^0 -1.5 + \frac{1}{\pi} \\
m_4^0 + 2.85
\end{bmatrix},
\end{equation}
where $\bm{k}^b = [k_1^b, k_2^b, k_3^b, k_4^b]^\top = [m_1^0 - \frac{2}{\pi}, m_2^0 + \frac{1}{\pi}, m_3^0 + \frac{1}{\pi}, m_4^0]^\top$ is defined by imposing that $\bm{m}_+^{0,0.5}(t) = \bm{m}^{0.5}$.

\textbf{Case $0.5 < t < 1$:} 
Using (\ref{eq:thetaAEx1Case3}) and that $\sin(\pi t) > 0$, $\cos(\pi t) < 0$, we have that
\begin{equation}
\bm{m}^{0.5,1}(t) = \\
\begin{bmatrix}
1.3 t + \frac{\cos(\pi t)}{\pi} + \frac{\sin(\pi t)}{\pi} \\
- \frac{\cos(\pi t)}{\pi} - 4t \\
- \frac{\sin(\pi t)}{\pi} - 3t \\
5.7 t 
\end{bmatrix}
+
\begin{bmatrix}
k^c_1 \\
k^c_2 \\
k^c_3 \\
k^c_4
\end{bmatrix},
\end{equation}
where $\bm{k}^c = [k^c_1, k^c_2, k^c_3, k^c_4]^\top = [m_1^0 - \frac{3}{\pi}, m_2^0 + \frac{1}{\pi}, m_3^0 + \frac{2}{\pi}, m_4^0]^\top$ is defined by imposing that $\bm{m}_-^{0.5,1}(t) = \bm{m}^{0.5}$.

\textbf{Case $t = 1$:}
Using (\ref{eq:thetaAEx1Case1}) and that $\cos(\pi t) < 0$, it follows that 
\begin{equation}
\begin{split}
\bm{m}^1  & = 
\begin{bmatrix}
\dot{m}_{4,1}t + \frac{\sin(\pi t)}{\pi}\\
- \dot{m}_{2,3}t\\
- \frac{\sin(\pi t)}{\pi} + \dot{m}_{2,3}t - \dot{m}_{3,4}t\\
\dot{m}_{3,4}t - \dot{m}_{4,1}t
\end{bmatrix} 
_{t = 1} 
+
\begin{bmatrix}
k^d_1 \\
k^d_2 \\
k^d_3 \\
k^d_4
\end{bmatrix} \\
& =
\begin{bmatrix}
m_1^0 + 1.3 - \frac{4}{\pi}\\
m_2^0 - 4 + \frac{2}{\pi}\\
m_3^0 - 3 + \frac{2}{\pi}\\
m_4^0 + 5.7
\end{bmatrix},
\end{split}
\end{equation}
where $\bm{k}^d = [k^d_1, k^d_2, k^d_3, k^d_4]^\top = [m_1^0 - \frac{4}{\pi}, m_2^0 + \frac{2}{\pi}, m_3^0 + \frac{2}{\pi}, m_4^0]^\top$ is defined by imposing that $\bm{m}_+^{0.5,1}(t) = \bm{m}^{1}$.

\textbf{Case $1 < t < 1.5$:}
Using (\ref{eq:thetaAEx1Case3}) and that $\sin(\pi t) < 0$, $\cos(\pi t) < 0$, we have that
\begin{equation}
\bm{m}^{1,1.5}(t) = \\
\begin{bmatrix}
1.3 t - \frac{\cos(\pi t)}{\pi} + \frac{\sin(\pi t)}{\pi} \\
\frac{\cos(\pi t)}{\pi} - 4t \\
- \frac{\sin(\pi t)}{\pi} - 3t \\
5.7 t 
\end{bmatrix}
+
\begin{bmatrix}
k^e_1 \\
k^e_2 \\
k^e_3 \\
k^e_4
\end{bmatrix},
\end{equation}
where $\bm{k}^e = [k^e_1, k^e_2, k^e_3, k^e_4]^\top = [m_1^0 - \frac{5}{\pi}, m_2^0 + \frac{3}{\pi}, m_3^0 + \frac{2}{\pi}, m_4^0]^\top$ is defined by imposing that $\bm{m}_-^{1,1.5}(t) = \bm{m}^{1}$.

\textbf{Case $t = 1.5$:}
Using (\ref{eq:thetaAEx1Case2}) and that $\sin(\frac{3\pi}{2}) < 0$, it follows that
\begin{equation}
\bm{m}^{1.5} = 
\begin{bmatrix}
\dot{m}_{4,1}t - \frac{\cos(\pi t)}{\pi} \\
\frac{\cos(\pi t)}{\pi} - \dot{m}_{2,3}t \\
\dot{m}_{2,3}t - \dot{m}_{3,4}t \\
\dot{m}_{3,4}t - \dot{m}_{4,1}t
\end{bmatrix}_{t=1.5}
+
\begin{bmatrix}
k_1^f \\
k_2^f \\
k_3^f \\
k_4^f
\end{bmatrix}
=
\begin{bmatrix}
 m_1^0 + 1.9 - \frac{6}{\pi}\\
 m_2^0 - 6 + \frac{3}{\pi}\\
 m_3^0 - 4.5 + \frac{3}{\pi}\\
m_4^0 + 8.5
\end{bmatrix},
\end{equation}
where $\bm{k}^f = [k_1^f, k_2^f, k_3^f, k_4^f]^\top = [m_1^0 - \frac{6}{\pi}, m_2^0 + \frac{3}{\pi}, m_3^0 + \frac{3}{\pi}, m_4^0]^\top$ is defined by imposing that $\bm{m}_+^{1,1.5}(t) = \bm{m}^{1.5}$.

\textbf{Case $1.5 < t < 2$:}
Using (\ref{eq:thetaAEx1Case3}) and that $\sin(\pi t) < 0$, $\cos(\pi t) > 0$, we have that
\begin{equation}
\bm{m}^{1.5,2}(t) = \\
\begin{bmatrix}
1.3 t - \frac{\cos(\pi t)}{\pi} - \frac{\sin(\pi t)}{\pi} \\
\frac{\cos(\pi t)}{\pi} - 4t \\
\frac{\sin(\pi t)}{\pi} - 3t \\
5.7 t 
\end{bmatrix}
+
\begin{bmatrix}
k^g_1 \\
k^g_2 \\
k^g_3 \\
k^g_4
\end{bmatrix},
\end{equation}
where $\bm{k}^g = [k^g_1, k^g_2, k^g_3, k^g_4]^\top = [m_1^0 - \frac{7}{\pi}, m_2^0 + \frac{3}{\pi}, m_3^0 + \frac{4}{\pi}, m_4^0]^\top$ is defined by imposing that $\bm{m}_-^{1.5,2}(t) = \bm{m}^{1.5}$.

\textbf{Case $t = 2$:}
Using (\ref{eq:thetaAEx1Case1}) and that $\cos(\pi t) > 0$, it follows that 
\begin{equation}
\begin{split}
\bm{m}^2 & = 
\begin{bmatrix}
\dot{m}_{4,1}t - \frac{\sin(\pi t)}{\pi}\\
- \dot{m}_{2,3}t\\
\frac{\sin(\pi t)}{\pi} + \dot{m}_{2,3}t - \dot{m}_{3,4}t\\
\dot{m}_{3,4}t - \dot{m}_{4,1}t
\end{bmatrix} 
_{t = 2} 
+
\begin{bmatrix}
k^h_1 \\
k^h_2 \\
k^h_3 \\
k^h_4
\end{bmatrix} \\
& =
\begin{bmatrix}
m_1^0 + 2.6 - \frac{8}{\pi}\\
m_2^0 - 8 + \frac{4}{\pi}\\
m_3^0 - 6 + \frac{4}{\pi}\\
m_4^0 + 11.4
\end{bmatrix},
\end{split}
\end{equation}
where $\bm{k}^h = [k^h_1, k^h_2, k^h_3, k^h_4]^\top = [m_1^0 - \frac{8}{\pi}, m_2^0 + \frac{4}{\pi}, m_3^0 + \frac{4}{\pi}, m_4^0]^\top$ is defined by imposing that $\bm{m}_+^{1.5,2}(t) = \bm{m}^{2}$.

\begin{figure}
\begin{subfigure}[t]{0.5\textwidth}
\centering
\includegraphics[width=\textwidth]{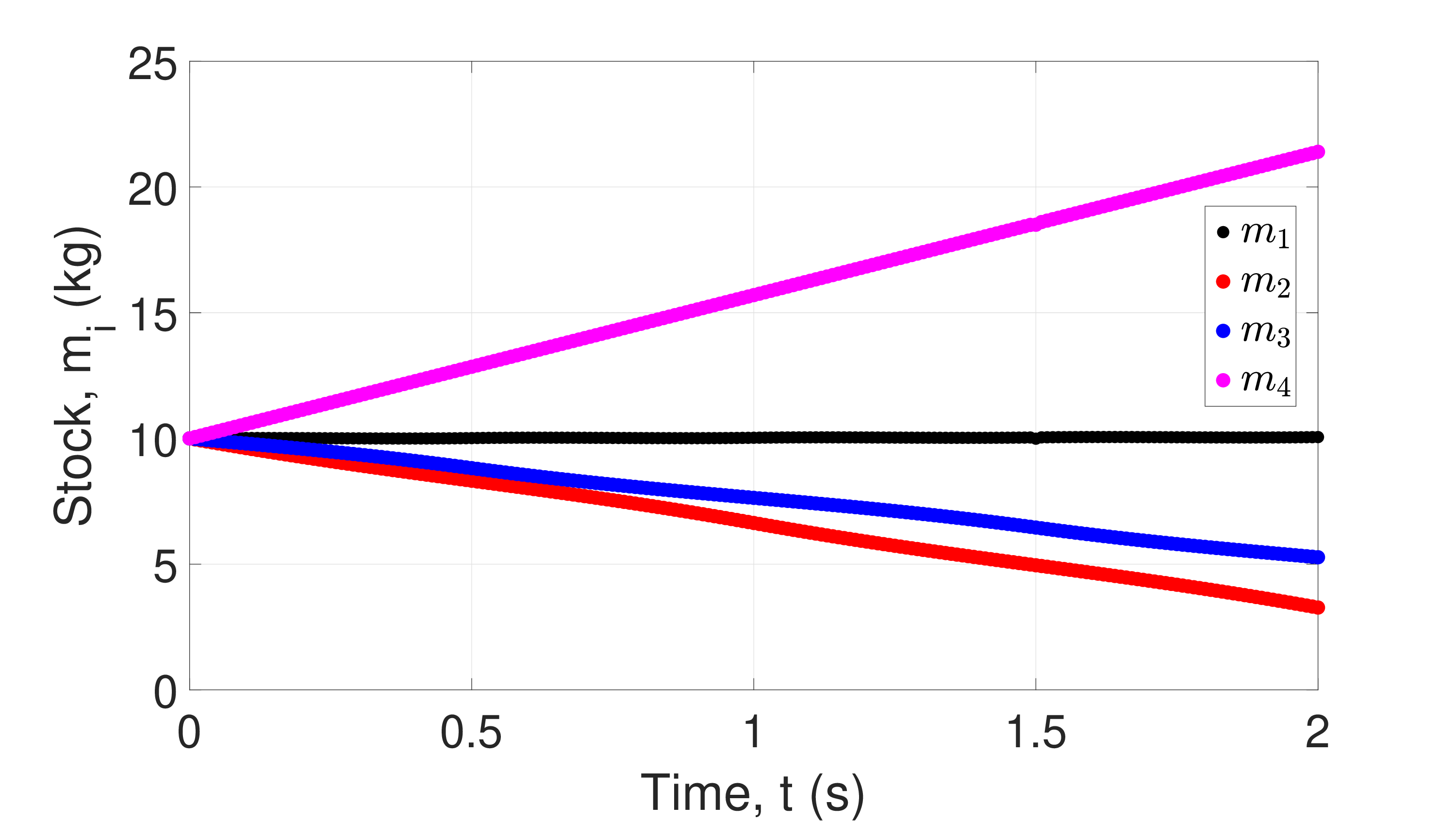}
\caption{Diagonal entries of $\bm{\Gamma}^{\text{mc}}(t)$}
\label{fig:Ex2-stocks}
\end{subfigure}
\begin{subfigure}[t]{0.5\textwidth}
\centering
\includegraphics[width=\textwidth]{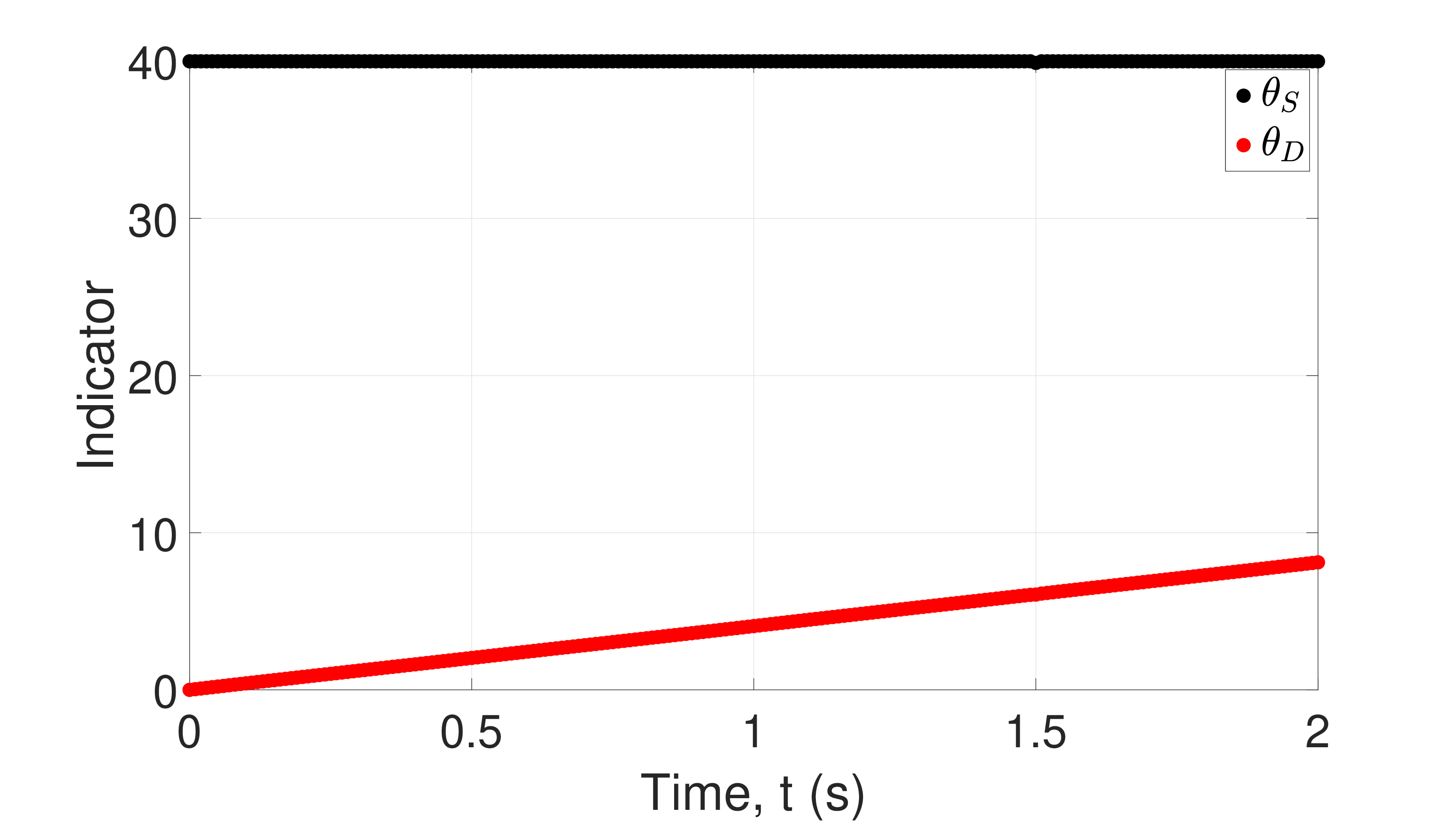}
\caption{$\theta_{\text{S}}$, $\theta_{\text{D}}$}
\label{fig:Ex2-indicators}
\end{subfigure}
\caption{Stocks and indicators in Example 2 for $[m_1^0, m_2^0, m_3^0, m_4^0]^\top = [10, 10, 10, 10]^\top$; the other indicators are unchanged and given in Example 1.}
\label{fig:Example2}
\end{figure}

The correct mass-flow matrix, which respects the mass conservation principle, is
\begin{equation}\label{eq:GammaExample2}
\bm{\Gamma}^{\text{mc}}(t) = 
\begin{bmatrix}
m^{\overline{0,2}}_1(t) & |\sin(\pi t)| & |\cos(\pi t)| & 0 \\
0 & m^{\overline{0,2}}_2(t) & 4 & 0 \\
0 & 0 & m^{\overline{0,2}}_3(t) & 7 \\
1.3 & 0 & 0 & m^{\overline{0,2}}_4(t)
\end{bmatrix},
\end{equation}
with $t \in [0, 2]$ and where $m_i^{\overline{t_a, t_b}}(t) = m_i(t)|_{{t_a} \leq t \leq {t_b}}$. The difference between the correct (\ref{eq:GammaExample2}) and the nonphysical (\ref{eq:GammaExample1}) mass-flow matrix is that the former has dynamic entries along the main diagonal, i.e., dynamic stocks, which follow from the imposition of the mass conservation principle (\ref{eq:constraintOfStocks}). The dynamic stocks are shown in Fig. \ref{fig:Ex2-stocks} for $\bm{m}^0 = [m_1^0, m_2^0, m_3^0, m_4^0]^\top = [10, 10, 10, 10]^\top$: the mass is equally distributed for $t = 0$, then it leaves $c^2_{2,2}$ and $c^3_{3,3}$ to accumulate in $c^4_{4,4}$ when $t > 0$; in contrast, the stock in $c^1_{1,1}$ remains constant. 

Following the correction of the mass-flow matrix, we need to update the indicators that depend on the diagonal entries (i.e., the mass stocks): $\theta_{\text{S}}(\bm{\Gamma})$ and $\theta_{\text{D}}(\bm{\Gamma})$. The two indicators are shown in Fig. \ref{fig:Ex2-indicators}. The total stock $\theta_{\text{S}}(\bm{\Gamma})$ is time-invariant since the network is a closed system, whereas the stock distribution $\theta_{\text{D}}(\bm{\Gamma})$ varies with time in contrast with Example 1 (see Fig. \ref{fig:Example1}), in which the mass conservation principle was violated. Specifically, the stock distribution is zero for $t = 0$ since $m_1^0 = m_2^0 = m_3^0 = m_4^0$, and then, it increases with time in agreement with the material accumulation in $c^4_{4,4}$ shown in Fig. \ref{fig:Ex2-stocks} for $t > 0$.       

\revision{
\section{Example for Solid Materials}\label{sec:caseStudy}
\subsection{Network Modeling and Simulation}
This section applies the TMN methodology \cite{zocco2023thermodynamical} to the case of solid plastics and measures its circularity using the proposed indicators. The considered system is a small network of single-use plastic for food and catering; the plastics of interest are polyethylene terephthalate (PET), high-density polyethylene (HDPE), and polypropylene (PP). An organization (indicated with $c^1_{1,1}$ in the compartmental digraph in Fig. \ref{fig:CaseStudy-graphs}) places an order with the manufacturer ($c^3_{3,3}$) of several products containing PET, HDPE, and PP; hence, the products are delivered to the organization by truck ($c^4_{3,1}$); after use, the plastic waste is collected from the organization and delivered to a recycling facility ($c^2_{2,2}$) by truck ($c^5_{1,2}$); finally, a fraction of the plastic is recycled and sent back to the organization by truck ($c^6_{2,1}$); factors such as contamination and material degradation result in non-recyclable plastic; these losses are considered in our model as it will be shown. The whole system can be depicted by the digraphs in Fig. \ref{fig:CaseStudy-graphs} having $n_v$ = 3, $n_a$ = 3, $n_c$ = 6, and corresponding to the network
\begin{equation}
\mathcal{N} = \{c^1_{1,1}, c^2_{2,2}, c^3_{3,3}, c^4_{3,1}, c^5_{1,2}, c^6_{2,1}\}.
\end{equation} 
\begin{figure}
\centering
\includegraphics[width=0.25\textwidth]{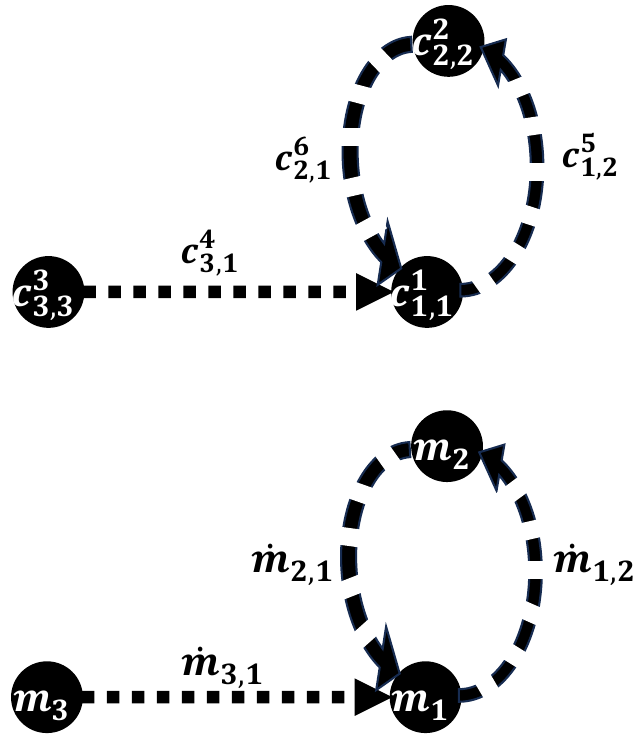}
\caption{\revision{Compartmental digraph (top) and mass-flow digraph (bottom) of the network considered in the example for solid plastics. Dashed arrows are for time-dependent flows.}}
\label{fig:CaseStudy-graphs}
\end{figure}
Note that, while several parameters of the model are taken from a real scenario, the life-cycle stages and their temporal sequence are a simplified representation of the reality. This study illustrates the use of the methodology and its difference from MFA: while MFA mainly relies on data analysis, TMNs mainly relies on dynamical systems; the data involved with TMNs enter the model as equation parameters.

We now need to model and simulate the dynamics of the flows and the stocks in the system respecting the mass balance principle. \emph{Since the considered plastics are solids transported in batches from a facility to another, the dynamics of stocks and flows is well described by a combination of discrete- and continuous-time terms. This is in contrast with the modeling of fluids transported in pipelines, which are well described by continuity and differential equations.} Hence, the dynamics of mass and flows in the network are depicted by the following hybrid equations, that is, equations involving both continuous- and discrete-time terms. The model parameters are summarized in Table \ref{tab:valuesOfParam}.            
\begin{equation}\label{eq:cs-m1}
\begin{gathered}
m_1(n+1) = m_1(n) + \overline{m}_{3,1} \delta_{t_{1,\text{in},4}}(t) \\ + \overline{m}_{2,1}\delta_{t_{1,\text{in},6}}(t) - \overline{m}_{1,2}\delta_{t_{1,\text{out},5}}(t), \quad n \in \overline{\mathbb{Z}}_+, 
\end{gathered}
\end{equation}
\begin{equation}
\begin{gathered}
m_2(n+1) = m_2(n) + \overline{m}_{1,2}\delta_{t_{2,\text{in},5}}(t) \\ - \overline{m}_{2,1}\delta_{t_{2,\text{out},6}}(t), \quad n \in \overline{\mathbb{Z}}_+,
\end{gathered}
\end{equation}
\begin{equation}\label{eq:cs-m3}
m_3(n+1) = m_3(n) - \overline{m}_{3,1}\delta_{t_{3,\text{out},4}}(t), \quad n \in \overline{\mathbb{Z}}_+,  
\end{equation}
\begin{equation}
\begin{gathered}
\dot{m}_{1,2}(n+1) = \frac{1}{T_{1,2}} m_{1,2}(n+1) \\ = \frac{1}{T_{1,2}}\left[m_{1,2}(n) + \overline{m}_{1,2} \left(\delta_{t_{5,\text{in},1}}(t) \right.\right. \\ \left.\left. - \delta_{t_{5,\text{out},2}}(t)\right)\right], \quad n \in \overline{\mathbb{Z}}_+,  
\end{gathered}
\end{equation}
\begin{equation}
\begin{gathered}
\dot{m}_{2,1}(n+1) = \frac{1}{T_{2,1}} m_{2,1}(n+1) \\ = \frac{1}{T_{2,1}}\left[m_{2,1}(n) + \overline{m}_{2,1} \left(\delta_{t_{6,\text{in},2}}(t) \right.\right. \\ \left.\left. - \delta_{t_{6,\text{out},1}}(t)\right)\right], \quad n \in \overline{\mathbb{Z}}_+,
\end{gathered}
\end{equation}
\begin{equation}\label{eq:cs-m31}
\begin{gathered}
\dot{m}_{3,1}(n+1) = \frac{1}{T_{3,1}} m_{3,1}(n+1) \\ = \frac{1}{T_{3,1}}\left[m_{3,1}(n) + \overline{m}_{3,1} \left(\delta_{t_{4,\text{in},3}}(t) \right.\right. \\ \left.\left. - \delta_{t_{4,\text{out},1}}(t)\right)\right], \quad n \in \overline{\mathbb{Z}}_+,
\end{gathered}
\end{equation}
\begin{equation}
\overline{m}_{2,1} = (1-x)\overline{m}_{1,2},
\end{equation}
where $n$ is the sample index, $\overline{\mathbb{Z}}_+$ is the set of nonnegative integers, and $\delta_{t_*}(t)$ is a rectangular pulse of short duration $\varepsilon$ emulating a unit impulse centered in $t_*$ and defined as
\begin{equation}\label{eq:deltaoft}
\delta_{t_*}(t) = \text{rect}\left(\frac{t - t_*}{\varepsilon}\right),
\end{equation}
where $\text{rect}(\cdot)$ is the rectangular function defined as \cite{vitettaTdS}
\begin{equation}
\text{rect}(\sigma) = 
\begin{cases}
1, & |\sigma| < 1/2 \\
1/2, & |\sigma| = 1/2 \\
0, & \text{otherwise}. \\
\end{cases}
\end{equation}
Moreover, $t_{k,\text{in},i}$ and $t_{k,\text{out},i}$ are the time instants at which the material enters the $k$-th compartment from the compartment $i$ and it leaves the $k$-th compartment for the compartment $i$, respectively, $\overline{m}_{i,j}$ is the mass of a batch of material transported from the compartment $i$ to the compartment $j$, $x \in [0, 1]$ is the fraction of material that is not recyclable due to contamination, and $T_{i,j}$ is the transportation time between the compartments $i$ and $j$.
\begin{table*}
\centering
\caption{\revision{Parameters used for the numerical study. The simulation starts at $t$ = 0 min and $n$ = 0. The results are shown in Fig. \ref{fig:csSimResults}. Note that these values are chosen to demonstrate the methodology and may not be realistic.}}
\label{tab:valuesOfParam}
\begin{tabular}{cccc} 
\hline
Compartment & Parameter & Value & Description\\ 
\hline
\multirow{7}{*}{$c^1_{1,1}$} & $m_{1,\text{PET}}(0)$ & 8.1 kg & Initial mass of PET in organization\\
& $m_{1,\text{HDPE}}(0)$ & 4.1 kg & Initial mass of HDPE in organization\\
& $m_{1,\text{PP}}(0)$ & 4.4 kg & Initial mass of PP in organization\\
& $t_{1,\text{out},5}$ & 100 min (1 h 40 min) & When plastics leave by truck for recycling\\
& $\overline{m}_{1,2,\text{PET}}$ & 8.1 kg & Mass of the batch of exiting PET \\ 
& $\overline{m}_{1,2,\text{HDPE}}$ & 4.1 kg & Mass of the batch of exiting HDPE \\
& $\overline{m}_{1,2,\text{PP}}$ & 4.4 kg & Mass of the batch of exiting PP \\  
\hline 
\multirow{5}{*}{$c^5_{1,2}$} & $T_{1,2}$ & 120 min (2 h 0 min) & Transportation time \\
& $t_{5,\text{out},2}$ & 220 min (3 h 40 min) & When plastics are unloaded from the truck and enter the recycling facility\\
& $m_{1,2,\text{PET}}(0)$ & 0 kg & Initial mass of PET in the truck\\
& $m_{1,2,\text{HDPE}}(0)$ & 0 kg & Initial mass of HDPE in the truck\\
& $m_{1,2,\text{PP}}(0)$ & 0 kg & Initial mass of PP in the truck\\
\hline 
\multirow{5}{*}{$c^2_{2,2}$} & $x$ & 0.3 & Fraction of non-recyclable plastics due to contamination \\
& $m_{2,\text{PET}}(0)$ & 20 kg & Initial mass of PET in recycling center\\
& $m_{2,\text{HDPE}}(0)$ & 20 kg & Initial mass of HDPE in recycling center\\
& $m_{2,\text{PP}}(0)$ & 20 kg & Initial mass of PP in recycling center\\
& $t_{2,\text{out},6}$ & 260 min (4 h 20 min) & When recycled plastics leave by truck for the organization\\
\hline 
\multirow{5}{*}{$c^6_{2,1}$} & $T_{2,1}$ & 30 min (0 h 30 min) & Transportation time \\
& $t_{6,\text{out},1}$ & 290 min (4 h 50 min) & When recycled plastics are unloaded from the truck and enter the organization\\
& $m_{2,1,\text{PET}}(0)$ & 0 kg & Initial mass of PET in the truck\\
& $m_{2,1,\text{HDPE}}(0)$ & 0 kg & Initial mass of HDPE in the truck\\
& $m_{2,1,\text{PP}}(0)$ & 0 kg & Initial mass of PP in the truck\\
\hline 
\multirow{7}{*}{$c^3_{3,3}$} & $m_{3,\text{PET}}(0)$ & 150 kg & Initial mass of PET in the manufacturer\\
& $m_{3,\text{HDPE}}(0)$ & 100 kg & Initial mass of HDPE in the manufacturer\\
& $m_{3,\text{PP}}(0)$ & 120 kg & Initial mass of PP in the manufacturer\\
& $t_{3,\text{out},4}$ & 1 min (0 h 1 min) & When new plastic products leave by truck for organization\\
& $\overline{m}_{3,1,\text{PET}}$ & 2 kg & Fraction of PET in new products ordered by organization\\ 
& $\overline{m}_{3,1,\text{HDPE}}$ & 0.5 kg & Fraction of HDPE in new products ordered by organization\\
& $\overline{m}_{3,1,\text{PP}}$ & 1 kg & Fraction of PP in new products ordered by organization\\
\hline 
\multirow{5}{*}{$c^4_{3,1}$} & $T_{3,1}$ & 6 min (0 h 6 min) & Transportation time \\
& $m_{3,1,\text{PET}}(0)$ & 0 kg & Initial mass of PET in the truck\\
& $m_{3,1,\text{HDPE}}(0)$ & 0 kg & Initial mass of HDPE in the truck\\
& $m_{3,1,\text{PP}}(0)$ & 0 kg & Initial mass of PP in the truck\\
&$t_{4,\text{out},1}$ & 7 min (0 h 7 min) & When new products are unloaded from the truck and enter the organization\\
\hline 
\end{tabular}
\end{table*}
Equations (\ref{eq:cs-m1})-(\ref{eq:cs-m31}) consider the total masses and flows, which result from the sum of the fractions of PET, HDPE, and PP, that is,
\begin{equation}
\begin{gathered}
m_i(n) = m_{i,\text{PET}}(n) + m_{i,\text{HDPE}}(n) \\ + m_{i,\text{PP}}(n), \quad i \in \{1,2,3\}, \, n \in \overline{\mathbb{Z}}_+,
\end{gathered}
\end{equation}
\begin{equation}
\begin{gathered}
m_{i,j}(n) = m_{i,j,\text{PET}}(n) + m_{i,j,\text{HDPE}}(n) \\ + m_{i,j,\text{PP}}(n), \quad i,j \in \{1,2,3\}, \, i \neq j, \, n \in \overline{\mathbb{Z}}_+,
\end{gathered}
\end{equation}
and
\begin{equation}
\begin{gathered}
\overline{m}_{i,j} = \overline{m}_{i,j,\text{PET}} + \overline{m}_{i,j,\text{HDPE}} \\ + \overline{m}_{i,j,\text{PP}}, \quad i,j \in \{1,2,3\}, \, i \neq j.
\end{gathered}
\end{equation}

Now, note that the exit of a batch of plastic from a compartment corresponds with the entrance of the same batch to the next compartment in the temporal sequence; meaning that, for example, the time in which a batch of PET leaves the manufacturer compartment corresponds to the time in which the same batch is loaded into the truck compartment leaving the manufacturer to deliver the batch. This leads to the imposition of the following equalities between the time instants in equations (\ref{eq:cs-m1})-(\ref{eq:cs-m31}):
\begin{equation}\label{eq:cs-constraints}
t_{i,\text{in},j} = t_{j,\text{out},i}, \quad i, j \in \{1, 2, 3, 4, 5, 6\}, \, i \neq j.
\end{equation}
Thus, in this case, the mass-flow matrix takes the form
\begin{equation}\label{eq:gammaMatCaseStudy}
\begin{gathered}
\bm{\Gamma}(\mathcal{N}) =
\begin{bmatrix}
\gamma_{1,1} & \gamma_{1,2} & \gamma_{1,3} \\
\gamma_{2,1} & \gamma_{2,2} & \gamma_{2,3} \\
\gamma_{3,1} & \gamma_{3,2} & \gamma_{3,3}  
\end{bmatrix}
\\ = 
\begin{bmatrix}
m_1(n+1) & \dot{m}_{1,2}(n+1) & 0 \\
\dot{m}_{2,1}(n+1) & m_2(n+1) & 0\\
\dot{m}_{3,1}(n+1) & 0 & m_3(n+1) 
\end{bmatrix}, \quad n \in \overline{\mathbb{Z}}_+,
\end{gathered} 
\end{equation} 
whose entries are given in (\ref{eq:cs-m1})-(\ref{eq:cs-m31}). A simulation of (\ref{eq:cs-m1})-(\ref{eq:cs-constraints}) is provided in Fig. \ref{fig:csSimResults} considering the parameters in Table \ref{tab:valuesOfParam}. The model was implemented in MathWorks Simulink 2024a using a variable-step discrete solver with a maximum step size of 0.4, with $\varepsilon$ = 0.4 min (see Equation (\ref{eq:deltaoft})), and it is publicly available\footnotemark[\value{footnote}].
\begin{figure*}
\begin{subfigure}{0.5\textwidth}
\centering
\includegraphics[width=\textwidth]{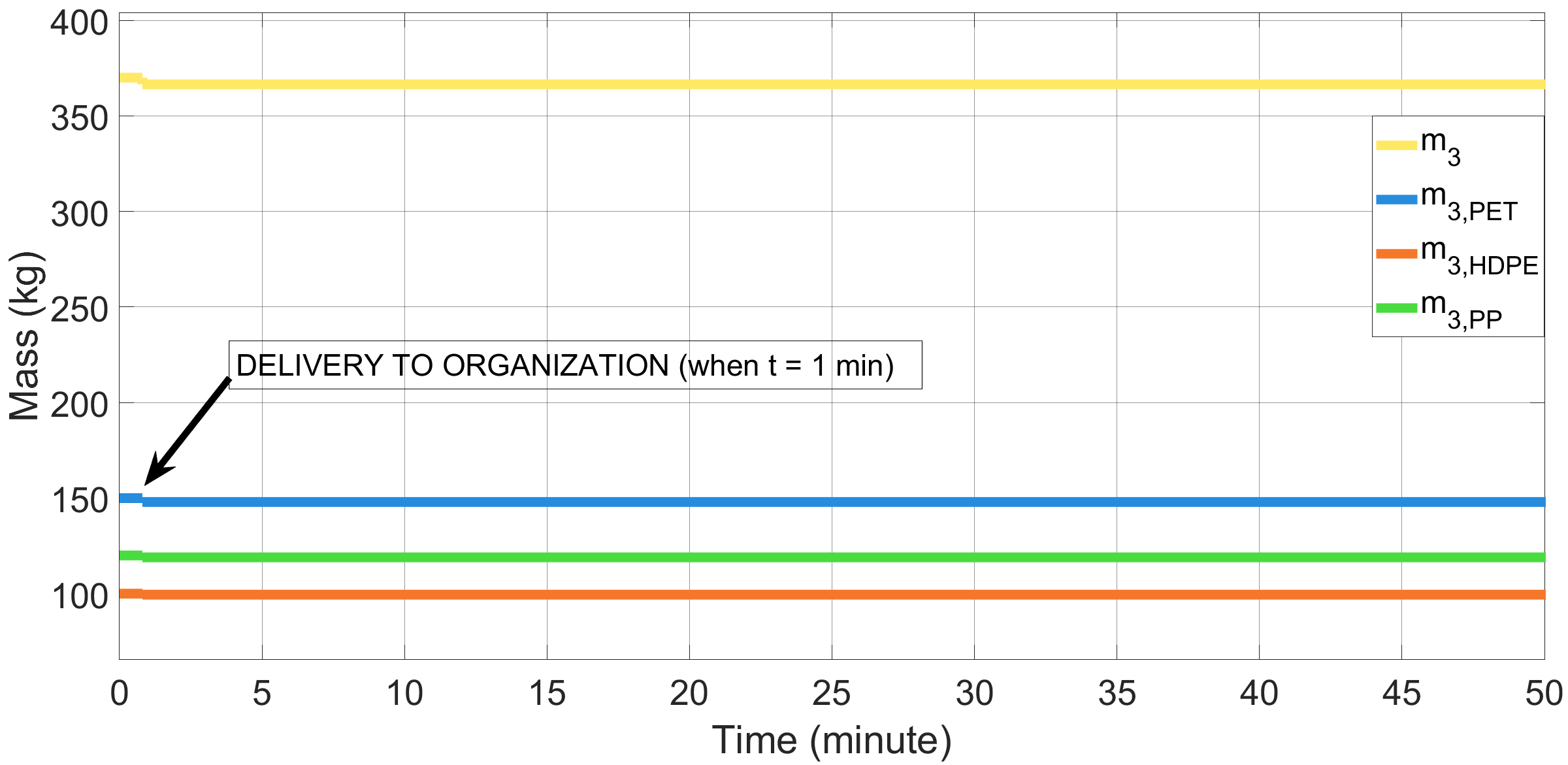}
\caption{Mass vs. time in $c^3_{3,3}$ (manufacturer)}
\label{fig:csManufacturer}
\end{subfigure}
\begin{subfigure}{0.5\textwidth}
\centering
\includegraphics[width=\textwidth]{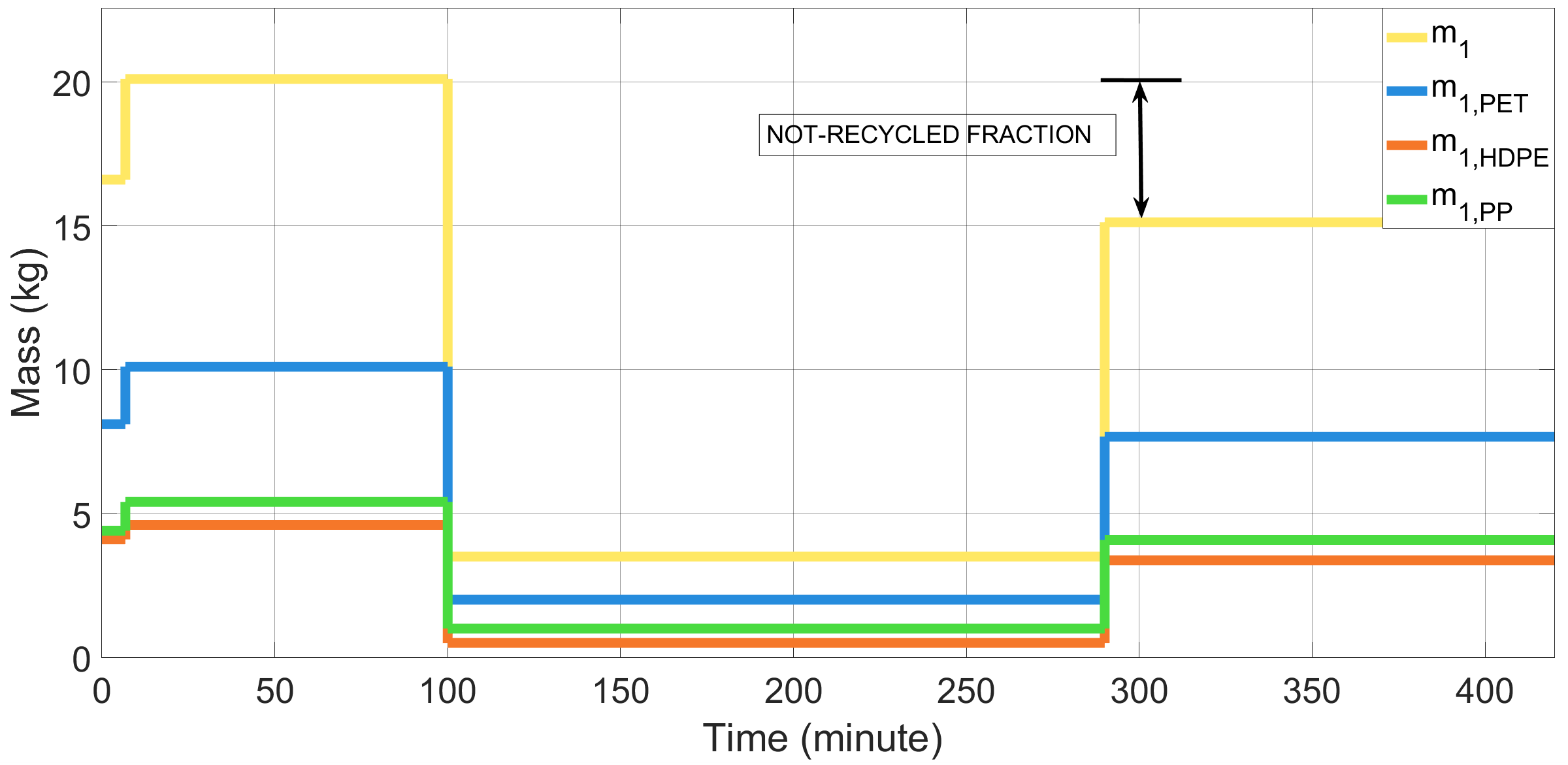}
\caption{Mass vs. time in $c^1_{1,1}$ (organization)}
\label{fig:csOrganization}
\end{subfigure}
\begin{subfigure}{0.5\textwidth}
\centering
\includegraphics[width=\textwidth]{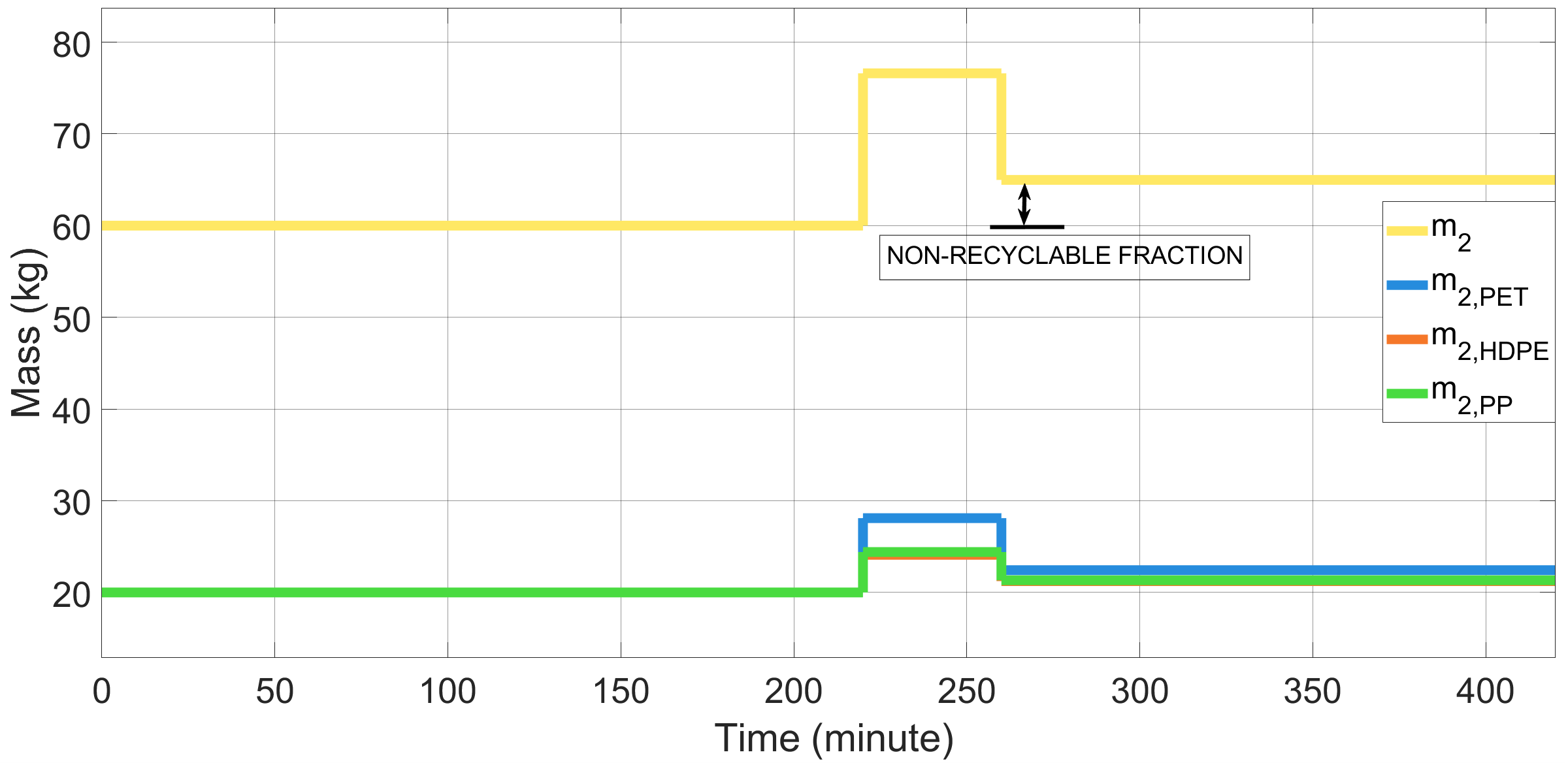}
\caption{Mass vs. time in $c^2_{2,2}$ (recycling facility)}
\label{fig:csRecycler}
\end{subfigure}
\begin{subfigure}{0.5\textwidth}
\centering
\includegraphics[width=\textwidth]{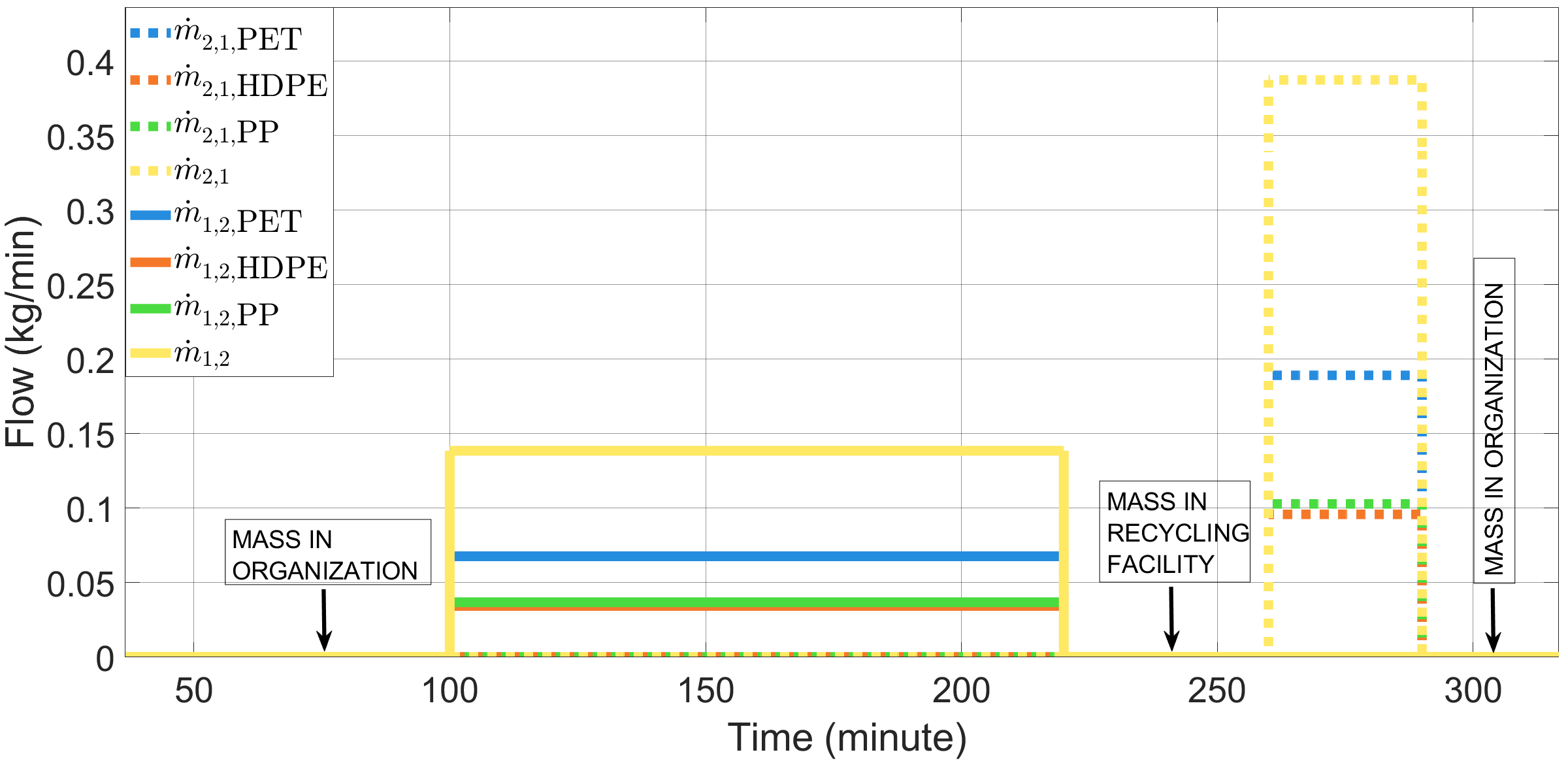}
\caption{Flow vs. time in $c^5_{1,2}$ and $c^6_{2,1}$ (transportations)}
\label{fig:csTransportations}
\end{subfigure}
\caption{\revision{Dynamics of stocks and flows in the example for solid materials. The model equations are given in (\ref{eq:cs-m1})-(\ref{eq:cs-constraints}), while the values of the parameters are in Table \ref{tab:valuesOfParam}. Yellow is for the total mass/flow, blue is for PET, orange is for HDPE, and green is for PP.}}
\label{fig:csSimResults}
\end{figure*}
As visible from Fig. \ref{fig:csOrganization}, the organization has an initial stock of PET, HDPE, and PP, which was accumulated over the previous two weeks (see for $t$ = 0 min). We assume that $t$ = 0 min corresponds to 10:00 am on the day at the end of the two-week period. At 10:00 am of that day, the organization places an order to the manufacturer (Fig. \ref{fig:csManufacturer}), which promptly delivers the new plastic products in 7 min since it is located near the organization. The delivery corresponds to the positive steps at $t_{4,\text{out},1}$ = 7 min in Fig. \ref{fig:csOrganization}. At $t_{1,\text{out},5}$ = 100 min, a truck collects some of the plastics in the organization. This event corresponds to the negative steps in Fig. \ref{fig:csOrganization} and the positive steps in Fig. \ref{fig:csTransportations}. Then, the truck transports the plastics to a recycling facility; the trip takes $T_{1,2}$ = 120 min due to a traffic congestion, which corresponds to the ``wide'' rectangles in Fig. \ref{fig:csTransportations}; in contrast, the ``tall'' rectangles in Fig. \ref{fig:csTransportations} correspond to the delivery of the recycled plastics back to the organization, whose transportation takes only $T_{2,1}$ = 30 min. During the time between the two groups of rectangles in Fig. \ref{fig:csTransportations}, the plastics are inside the recycling center. To keep the simulation horizon of limited length and capture all the dynamics within the figures, we assume that the recycling process takes only 40 min (Fig. \ref{fig:csRecycler}). Once the recycling is complete, the plastics are sent back to the organization by truck at $t_{2,\text{out},6}$ = 260 min, i.e., 4 h and 20 min after 10:00 am, which corresponds to $t = 0$. The delivery of the recycled plastics to the organization corresponds to the positive steps at $t_{1,\text{in},6}$ = 290 min in Fig. \ref{fig:csOrganization}. The fraction of non-recyclable plastics due to contamination (expressed with $x$ = 0.3) is indicated in Figs. \ref{fig:csOrganization} and \ref{fig:csRecycler}. Note that PET, HDPE, and PP are considered to enter/exit the compartments at the same time, thus their positive/negative steps overlap (see, for example, $t$ = 100 min in Fig. \ref{fig:csOrganization}).

\subsection{Indicators and Discussion} 
Several proposed indicators, namely, $\lambda_{\text{Y}}$, $\lambda_{\text{GS}}$, $\lambda_{\text{GT}}$, $\lambda_{\text{HS}}$, $\lambda_{\text{HT}}$, $\lambda_{\text{AS}}$, $\lambda_{\text{AT}}$, and $\lambda_{\text{S}}$, are highly affected by the \emph{simultaneous existence of non-null flows within cycles} such as those in the examples covered in Section \ref{sec:Examples}. This network, instead, does not show flows that are simultaneously non-null as visible in Fig. \ref{fig:csTransportations}, where the two flows $\dot{m}_{1,2}$ and $\dot{m}_{2,1}$ do not overlap in time. In other words, in this network, there are no cycles, i.e.,
\begin{equation}\label{eq:noLoopsInCS}
n_\phi(t) = 0,\quad \forall t,     
\end{equation}  
despite the fact that the digraph in Fig. \ref{fig:CaseStudy-graphs} seems to have one loop involving $c^1_{1,1}$ and $c^2_{2,2}$; since that loop has time-dependent flows which are non-null at different times, Equation (\ref{eq:noLoopsInCS}) holds.
As a consequence, we have that $\lambda_{\text{Y}}$ = 0, and hence, $\lambda_{\text{GS}}$ = $\lambda_{\text{GT}}$ = $\lambda_{\text{HS}}$ = $\lambda_{\text{HT}}$ = $\lambda_{\text{AS}}$ = $\lambda_{\text{AT}}$ = $\lambda_{\text{S}}$ = 0. For the average connectivity we have
\begin{equation}\label{eq:csLambdaC}
\begin{gathered}
\lambda_{\text{C}}(t) = 
\begin{cases}
\frac{2}{3} \quad \text{if } t_{3,\text{out},4} \leq t \leq t_{4,\text{out},1}, \\ \quad\,\,\,\, t_{1,\text{out},5} \leq t \leq t_{5,\text{out},2}, \\ \quad\,\,\,\, t_{2,\text{out},6} \leq t \leq t_{6,\text{out},1},\\
0 \quad \text{otherwise},\\
\end{cases}
\end{gathered}
\end{equation} 
while in $\lambda_{\text{D}}$ occurs the division by zero, and specifically,    
\begin{equation}\label{eq:csLambdaD}
\begin{gathered}
\lambda_{\text{D}}(t) = 
\begin{cases}
\frac{0}{0} \quad \text{if } t < t_{3,\text{out},4}, \\ \quad\,\,\,\, t_{4,\text{out},1} < t < t_{1,\text{out},5}, \\ \quad\,\,\,\, t_{5,\text{out},2} < t < t_{2,\text{out},6}, \\ \quad\,\,\,\, t > t_{6,\text{out},1}, \\
0 \quad \text{if } t_{3,\text{out},4} \leq t \leq t_{4,\text{out},1}, \\ \quad\,\,\,\, t_{2,\text{out},6} \leq t \leq t_{6,\text{out},1},\\
\infty \quad \text{otherwise}.\\
\end{cases}
\end{gathered}
\end{equation} 
The auxiliary indicators $\theta_{\text{S}}$, $\theta_{\text{F}}$, and $\theta_{\text{D}}$ are depicted in Fig. \ref{fig:csIndicators}, while $\bm{\theta}_{\text{A}}$ cannot be calculated in this case of solids since the dynamical equations of the stocks (\ref{eq:cs-m1})-(\ref{eq:cs-m3}) involve both continuous- and discrete-time terms, and hence, their derivative with respect to $t$ is not defined. 
\begin{figure*}
\begin{subfigure}{0.33\textwidth}
\centering
\includegraphics[width=\textwidth]{Figures/theta\_S}
\caption{$\theta_\text{S}$}
\label{fig:thetaS}
\end{subfigure}
\begin{subfigure}{0.33\textwidth}
\centering
\includegraphics[width=\textwidth]{Figures/theta\_F}
\caption{$\theta_\text{F}$}
\label{fig:thetaF}
\end{subfigure}
\begin{subfigure}{0.33\textwidth}
\centering
\includegraphics[width=\textwidth]{Figures/theta\_D}
\caption{$\theta_\text{D}$}
\label{fig:thetaD}
\end{subfigure}
\caption{\revision{Auxiliary indicators for the example of solid materials.}}
\label{fig:csIndicators}
\end{figure*}
This numerical example highlights a characteristic of the formulation of the indicators $\lambda_{\text{Y}}$, $\lambda_{\text{GS}}$, $\lambda_{\text{GT}}$, $\lambda_{\text{HS}}$, $\lambda_{\text{HT}}$, $\lambda_{\text{AS}}$, $\lambda_{\text{AT}}$, and $\lambda_{\text{S}}$, which was not visible with the examples in Section \ref{sec:Examples}; the characteristic is that these indicators are equal to zero unless there are \emph{simultaneous flows closed in a loop}. In practice, this happens if batches of material are transported simultaneously in a loop, i.e., 
\begin{equation}\label{eq:conditionForLoopInCS}
\exists t : \dot{m}_{1,2}(t) \neq 0 \land \dot{m}_{2,1}(t) \neq 0. 
\end{equation}
With continuous-time flows typical of water networks, it is common to find some $t$ such that (\ref{eq:conditionForLoopInCS}) holds. In contrast, the transportation of solids as in this case is carried out in batches, and hence, the flows are non-null only within the time window during which the transportation takes place; this, in general, makes the satisfaction of (\ref{eq:conditionForLoopInCS}) harder. This phenomenon affects also $\lambda_{\text{C}}$ (\ref{eq:csLambdaC}), which differs from zero only when a transportation is performed. The directionality $\lambda_{\text{D}}$ (\ref{eq:csLambdaD}) yields 0/0 when there is no transportation, whereas $\lambda_{\text{D}} = \infty$ when $t_{1,\text{out},5} \leq t \leq t_{5,\text{out},2}$, while it equals zero anywhere else. The total stock $\theta_{\text{S}}$ (Fig. \ref{fig:thetaS}) is 446.6 kg for $t$ = 0 (the sum of the initial stock in $c^1_{1,1}$, $c^2_{2,2}$, and $c^3_{3,3}$) and at any times in which there is no transportation; in correspondence of the three transportations, $\theta_{\text{S}} < 446.6$ kg because part of the mass is on the trucks. In contrast, the total flow $\theta_{\text{F}}$ (Fig. \ref{fig:thetaF}) is different from zero only during the transportation stages (the three rectangular pulses). Finally, the stock distribution $\theta_{\text{D}}$ (Fig. \ref{fig:thetaD}) reaches the maximum of 195 kg when $\dot{m}_{1,2} \neq 0$ and the minimum of 190 kg when $t_{4,\text{out},1} < t < t_{1,\text{out},5}$; variations of $\theta_{\text{D}}$ occur in proximity of every exchange of mass between compartments, i.e., when $t \approx 1, 7, 100, 220, 260, 290$ min.            

The fact that several indicators equal zero, namely, $\lambda_{\text{GS}}$ = $\lambda_{\text{GT}}$ = $\lambda_{\text{HS}}$ = $\lambda_{\text{HT}}$ = $\lambda_{\text{AS}}$ = $\lambda_{\text{AT}}$ = $\lambda_{\text{S}}$, and $\lambda_{\text{Y}}$, indicates that the circularity of the network should be improved. In order to increase it, as pointed out before, simultaneous flows within cycles must occur. In practice, this means that two or more trucks must carry out their delivery at the same time. This can happen if transportations are more frequent than the one modeled in this case, which is common in real scenarios. Indeed, in reality, the delivery of products and the collection of waste for recycling is performed by multiple trucks simultaneously covering different routes, whereas in this numerical example there is only one truck in motion at a time. Using more realistic model parameters and increasing the size of the network is left as future work; this example has provided the modeling principles to accurately depict stocks and flows with dynamics faster than 1 minute.

\section{Conclusion}\label{sec:Concl}
This paper has proposed a graph-based measure of material flow circularity leveraging the formalism of thermodynamical material networks to improve the theoretical foundations of circular flow modeling and design. Our illustrative examples covered the calculation of the indicators for the case of dynamic continuous-time flows and provided simulations coherent with the network architecture. The diagonal and off-diagonal entries of the mass-flow matrix are constrained by the mass conservation principle, whose satisfaction has been the critical step. If the principle is violated as in the first example, the circularity indicators are nonphysical. 

Subsequently, we considered a numerical example for solid plastics. Modeling the system as a TMN has shown to capture dynamics of stocks and flows faster than 1 minute such as the entrance of a batch of material into the recycling center or the departure of a truck carrying out a delivery. The model equations required a total of 34 parameters (the ones in Table \ref{tab:valuesOfParam}), which is less than the amount of data required to achieve such a modeling accuracy using MFA. The reduction in the amount of data is achieved by carefully looking at the actual behavior of stocks and flows in reality and replicate it using hybrid dynamical equations, i.e., equations containing both continuous- and discrete-time terms. In contrast, MFA would require \emph{the full record of historical data} to detect the events of entry/exit of material depicted in this example if they were real. Specifically, MFA would require the following number of data: the events in this example occur at specific minutes, and hence, a sampling time of 1 minute is the maximum sampling time that guarantees the detection of all the events (meaning that a sample is recorded every minute); the first event occurs at $t = t_{3,\text{out},4} = 1$ min, while the last one occurs at $t = t_{6,\text{out},1} = 290$ min; thus, all the events to be detected take place within the time window  $\Delta = t_{6,\text{out},1} - t_{3,\text{out},4} = 289$ min; with a sampling time of 1 min, this yields that 289 samples are needed to record a stock/flow for not missing any entry-exit event; now, how many stocks/flows need to be recorded with MFA to plot all the events as we did here? It is needed the recording of all the stocks for each material type, while the flows can be derived from mass balances; since there are 3 stocks (one in the manufacturer, one in the organization, and one in the recycler) per each of the 3 types of plastic, the total number of data needed with MFA to detect all the entry-exit events is $3 \times 3 \times 289 = 2,601$ historical data, which is 2,601/34 = 76 times the amount of data we needed to model all the events with our approach (we needed only the 34 parameters in Table \ref{tab:valuesOfParam}).        

Once the network was modeled, the circularity indicators were calculated; most of the indicators equal zero because this example for solid materials has no simultaneous flows within cycles. This occurs because the example is a simplified version of real networks, and yet, it has required more than 12 equations to be defined. This example provides the principles for capturing and simulating the real dynamics of stocks and flows. Scaling-up the model is left as future work to provide networks useful to organizations and industry for the design of circular flows of materials. Other future work will be to formulate the indicators in such a way that they do not require simultaneous flows within cycles, but rather flows in cycles within a fixed time window, e.g., a day.

\section*{Acknowledgements}
Federico Zocco thanks Pantelis Sopasakis, Centre for Intelligent Autonomous Manufacturing Systems, School of Electronics, Electrical Engineering and Computer Science, Queen's University Belfast, UK, and Beatrice Smyth, Research Centre in Sustainable Energy, School of Mechanical and Aerospace Engineering, Queen's University Belfast, UK, for the valuable comments; he also thanks Bronagh Millar, Advanced Manufacturing Innovation Centre (AMIC), Queen's University Belfast, UK, for the inputs regarding supply-recovery chains of plastics. 
}

\ifCLASSOPTIONcaptionsoff
  \newpage
\fi

\bibliographystyle{IEEEtran}
\bibliography{References}

\end{document}